\documentclass[11pt]{article}

\usepackage{amssymb,amsmath}
\usepackage{graphicx,psfrag,color}
\usepackage{fullpage}
\usepackage{hyperref}

\begin{document}

\title{Comparing theories: the dynamics of changing
vocabulary.\\ A case-study in relativity theory.}
\author{Andr\'eka, Hajnal and N\'emeti, Istv\'an%
\footnote{Research supported by the Hungarian grant for basic
research OTKA No K81188.}\\
Alfr\'ed R\'enyi Institute of Mathematics, Hungarian Academy of
Sciences,\\ 1053 Budapest, Re\'altanoda u.\ 13-15, Hungary,\\
andreka.hajnal@renyi.mta.hu, nemeti.istvan@renyi.mta.hu}
 \maketitle

\newtheorem{thm}{Theorem}[section]
\newenvironment{theorem}
{\begin{thm}\it}{\end{thm}}

\definecolor{gray}{gray}{.6}
\newcommand{\grey}[1]{\textcolor{red}{#1}\/}
\newcommand{\draft}[1]{\textcolor{gray}{#1}}
\newcommand{\mar}[1]{\textcolor{red}{#1}}

\newcommand{\SigTh}{\mbox{\textsc{SigTh}}} %
\newcommand{\SpecRel}{\mbox{\textsc{SpecRel}}}
\newcommand{\EFd}{\mbox{\textsc{EFd}}}
\newcommand{\Compl}{\mbox{\textsc{Compl}}} %
\newcommand{\Q}{\mbox{\small \sf Q}} %
\newcommand{\B}{\mbox{\small \sf B}} %
\newcommand{\Obs}{\mbox{\small \sf Obs}} %
\newcommand{\Ph}{\mbox{\small \sf Ph}} %
\newcommand{\W}{\mbox{\small \sf W}} %
\newcommand{\AxPh}{\mbox{\small \sf AxPh}} %
\newcommand{\AxEv}{\mbox{\small \sf AxEv}} %
\newcommand{\AxSelf}{\mbox{\small \sf AxSelf}} %
\newcommand{\AxFd}{\mbox{\small \sf AxFd}} %
\newcommand{\AxSym}{\mbox{\small \sf AxSym}} %
\newcommand{\AxUp}{\mbox{\sf Ax$\uparrow$}} %

\newcommand{\Sig}{\mbox{\small \sf Sig}} 
\newcommand{\Par}{\mbox{\small \sf Par}} 
\newcommand{\Tran}{\mbox{\small \sf T}} 
\newcommand{\T}{\Tran} 
\newcommand{\Rec}{\mbox{\small \,\sf R\,}} 
\newcommand{\R}{\Rec} 
\newcommand{\Meet}{\mbox{\small \sf Meet}}
\newcommand{\meet}{\mbox{\small \sf meet}}

\newcommand{\F}{\mbox{\sf F}}
\newcommand{\Ev}{\mbox{\small \sf Ev}} %
\newcommand{\Space}{\mbox{\small \sf Space}} %
\newcommand{\Time}{\mbox{\small \sf Time}} %
\newcommand{\tim}{\mbox{\small \sf time}} %
\newcommand{\loc}{\mbox{\small \sf loc}} %
\newcommand{\cord}{\mbox{\small \sf cor}} %
\newcommand{\simu}{\equiv} %

\newcommand{\e}{\varepsilon} %
\newcommand{\f}{\varphi} %
\newcommand{\s}{\mbox{\small $\sigma$}} %
\newcommand{\g}{\mbox{\small $\gamma$}} %
\newcommand{\ex}{\mbox{e}} 
\newcommand{\so}{\mbox{\small $o$}}
\newcommand{\x}{\mbox{\small $\xi$}}

\newcommand{\dd}{\mbox{\small $\delta$}}
\newcommand{\uu}{\mbox{\small $\iota$}}
\newcommand{\oo}{\mbox{\small $o$}}
\newcommand{\st}{\mbox{\small $\tau$}}
\newcommand{\ea}{\mbox{$e$}}
\newcommand{\LL}{\mathcal L} %

\newcommand{\rep}{\mbox{\small \sf rep}}
\newcommand{\Tu}{\mbox{\small \sf Tu}}
\newcommand{\Beg}{\mbox{\small \sf Beg}} %
\newcommand{\End}{\mbox{\small \sf End}} %
\newcommand{\Ed}{\mbox{\small \sf Ed}} %
\newcommand{\Fp}{\mbox{\small \sf Fp}} %
\newcommand{\Op}{\mbox{\small \sf Op}} %
\newcommand{\Ted}{\mbox{\small \sf Edt}} %
\newcommand{\Bw}{\mbox{\small \sf Bw}} %
\newcommand{\M}{\mbox{\small \sf M}} %
\newcommand{\Col}{\mbox{\small \sf Col}} %
\renewcommand{\S}{\mbox{\small \sf S}} %
\newcommand{\Ort}{\mbox{\small \sf Ort}} %
\newcommand{\pa}{\mbox{\small \sf pa}} 
\newcommand{\de}{\ \ :=\ \ } %
\newcommand{\deiff}{{\ \ :\Leftrightarrow\ \ }} %
\newcommand{\Mm}{\mathfrak{M}}
\newcommand{\Cc}{\mathfrak{C}}
\newcommand{\sMm}{\mbox{\tiny $\mathfrak{M}$}}
\newcommand{\Mod}{\mbox{\small \sf Mod}} %
\newcommand{\Th}{\mbox{\small \sf Th}} %
\newcommand{\var}{\mbox{\small \sf var}} %
\newcommand{\Var}{\mbox{\small \sf Var}} %
\newcommand{\tr}{\mbox{\small \sf tr}} %
\newcommand{\Tr}{\mbox{\small \sf Tr}} %
\newcommand{\Points}{\mbox{\small \sf Points}} %
\newcommand{\Lines}{\mbox{\small \sf Lines}} %
\newcommand{\MS}{\mbox{\small \sf MS}} %
\newcommand{\is}{\mbox{\small \sf iso}} %
\newcommand{\Med}{\mbox{\small \sf Edr}} %
\newcommand{\Pm}{\mbox{\small \sf Pm}} %
\newcommand{\wl}{\mbox{\small \sf wl}} %
\newcommand{\wls}{\mbox{\small \sf wls}} %
\newcommand{\wla}{\mbox{\small \sf wlp}} %
\newcommand{\ev}{\mbox{\small \sf ev}} %
\newcommand{\AfP}{\mbox{\textsc{AfP}}} %
\newcommand{\uv}{\mbox{\small \sf Uv}} %
\newcommand{\eq}{\mbox{\small \sf Eq}} %
\newcommand{\ph}{\mbox{\scriptsize \sf Ph}} %
\newcommand{\obs}{\mbox{\scriptsize \sf Obs}} %
\newcommand{\pj}{\mbox{\small \sf pj}}


\begin{abstract}There are several first-order logic (FOL)
axiomatizations of special relativity theory in the literature, all
looking different but claiming to axiomatize the same physical
theory. In this paper, we elaborate a comparison, in the framework
of mathematical logic, between these FOL theories for special
relativity. For this comparison, we use a version of mathematical
definability theory in which new entities can also be defined
besides new relations over already available entities. In
particular, we build an interpretation (in Tarski's sense) of the
reference-frame oriented theory $\SpecRel$ of \cite{AMNHbSL} into
the observationally oriented Signalling theory of James
Ax~\cite{Ax}. This interpretation provides $\SpecRel$ with an
operational/experimental semantics. Then we make precise,
``quantitative" comparisons between these two theories via using the
notion of definitional equivalence. This is an application of logic
to the philosophy of science and physics in the spirit of van
Benthem's \cite{vB82, vB12}.\end{abstract}

\section{Introduction} \label{sec:intro}

This paper is about an application of logic to the methodology of
science in the spirit of van Benthem's \cite{vB82, vB12}.

There are several axiomatizations of special relativity theory
available in the literature, all looking different but claiming to
axiomatize the same physical theory. Such are, among many others,
the ones in Andr\'eka et al.~\cite{Vienna}, Ax~\cite{Ax},
Goldblatt~\cite{Gold}, Schelb~\cite{Schelb}, Schutz~\cite{Schutz},
Suppes~\cite{Sup59}. These papers talk about very different kinds of
objects: \cite{Vienna} talks about reference frames, \cite{Ax} talks
about particles and signals, \cite{Gold} looks like a purely
geometrical theory about orthogonality, the central notion of
\cite{Sup59} is the so-called Minkowski-metric, etc. While, as
usual, one gets a better picture of this area via a variety of
different ``eyeglasses", the following questions arise. What are the
connections between these theories? Do they all talk about the same
thing? If they do, do they capture it to the same extent, or is one
axiomatization more detailed or accurate than some of the others? In
this paper we want to show how, in the framework of mathematical
logic, a concrete, tangible comparison/connection can be elaborated
between these theories for special relativity and what we can gain
from such an investigation.

For this comparison, we have to use a form of logical definability
theory in which totally new kinds of entities can be defined as
opposed to traditional definability theory where only new relations
can be defined on already available entities. The existing methods
of definability theory had to be modified and refined for the
purposes of the present situation. Thus, definability theory, too,
profits from such an application. In the present paper, we elaborate
in detail one piece of such a comparison: we construct an
interpretation of the relativity theory in \cite{Vienna} talking
about reference frames into the theory in \cite{Ax} which talks
about particles emitting and absorbing signals. Then we construct an
inverse interpretation for showing which versions of the two
theories are definitionally equivalent. Since this is a case-study
showing applicability of the proposed method for
comparing/connecting theories, we tried to give all the detail
needed. This is why some sections of the paper may look somewhat
technical.

An insight of last century mathematical logic is that it is
important to fix the vocabulary of a first-order logic theory and
stay inside the so obtained language while in a specific theory
(see, e.g., \cite{Tarski36}). The symbols in the vocabulary%
\footnote{Other names for vocabulary are signature and set of
nonlogical constants} are the concepts that are not analyzed further
in the given theory, they are called thus basic (or primitive)
concepts. But this is not a forever frozen state: we may decide to
analyze further the basic concepts of this vocabulary and we can do
this in the form of building an interpretation (in the sense of
mathematical logic) into another language the vocabulary of which
consists of new basic concepts, and the interpretation gives us the
information of how the ``old basic concepts" are built up from the
``new basic concepts" as refined ones. The interpretation we
construct in this paper thus refines the basic concept of a
reference frame in terms of just sending and receiving signals. To
refine the basic concepts of this Signalling theory, we can
interpret it to, say, in a theory of electromagnetism, or in a
quantum-mechanical theory.

Such an interpretation may also be regarded as defining a so-called
operational semantics for the basic concepts of the first theory.
Starting with the Vienna Circle, several authors suggest that a
physical theory is a more complex object than just a set of
first-order logic (FOL for short) formulas. A physical theory, they
propose, is a FOL theory together with instructions for how to
interpret the basic symbols (or vocabulary) of this theory ``in the
real physical world". (Following Carnap~\cite{Carnap}, this is often
called a ``(partially) interpreted theory".) We want to show in
section \ref{sec:reduc} that such an ``operational semantics" can be
taken to be an interpretation in the sense of mathematical logic.

Returning to our concrete example, an operational semantics should
say something about how we obtain or set up (in the real world) the
reference frames for special relativity theory. Usually, rigid
meter-rods and standard clocks are used for this purpose (e.g.,
\cite{Wheeler}). However, as \cite{Szabo} points out, we cannot use
these rigid meter-rods in astronomy or cosmology. The interpretation
we give in detail in this paper results also in an
operational/experimental/observational definition for setting up a
reference frame by just relying on sending and receiving
light-signals. This method can be used, in principle, in the above
mentioned astronomical scale.

Summing up, the first language in an interpretation has the
theoretical concepts while the basic concepts of the second language
are the observational ones (for the observational-theoretical
duality see, e.g., \cite{vB82,Frie}). We can look at the same
interpretation ``from the other direction": In our example, we may
imagine someone living in a space-time, exploring his surroundings
by sending and receiving signals, and during this process, he
devises so-called theoretical concepts which make thinking more
efficient. In particular, he may devise the concept of a reference
frame, and even the concept of quantities forming a field, as mental
constructs having concrete definitions in terms of observations. The
tools of mathematical logic, and more closely those of definability
theory (interpretations are among them) can be used for modeling
this emergence of theoretical concepts.

A further aim of the present approach of comparing theories is
shifting the emphasis from working inside a single huge theory to
working in a modularized hierarchy of smaller theories connected in
many ways. Usually this approach is called theory-hierarchy. We note
that this is not so much a hierarchy as rather a category of
theories, technically the category of all FOL-theories as objects
with interpretations as morphisms of the category. This direction of
replacing a huge theory with a category of small theories is present
in many parts of science. In foundational thinking, \cite{FOM}
emphasizes this. In computer science, it is present in the form of
structured programming. ``Putting theories together" of Burstall and
Goguen \cite{BG} refers to the act of computing/generating colimits
of certain diagrams in this category. Even in such practical areas
as using a huge medical data-base the need of modularizing arises:
it is necessary to ``break up" the given data-base and generate many
smaller ones according to the query at hand \cite{Volt, Lutz}. The
interpretation going from special relativity as formalized in
\cite{Vienna} into the more observational Signalling theory of
\cite{Ax} we build in the present paper is but one morphism of this
huge dynamic category of FOL theories.

The content of this paper can also be viewed as preparing the ground
for an application of algebraic logic to relativity theory, as
follows. The cylindric algebra of a theory is an abstract
representation of the structure of concepts expressible/definable in
that theory and a homomorphism between two cylindric algebras
corresponds to an interpretation between the corresponding theories.
Hence the category of all FOL theories is basically the same as the
category of cylindric algebras as objects and homomorphisms as
morphisms. There are, for example, well known and understood methods
for how to compute colimits in this category of algebras.

There are still many questions and phenomena to be understood in
this area of application of logic. For example, what are the
desirable or good properties of an interpretation for being
informative about the theories in one or other respect? Consider for
a second the definability/interpretability picture between
scientific theories (in FOL) in two versions: (1) in the framework
of traditional definability theory, and (2) in the new, extended
theory of definability used in the present paper. What are the
characteristic differences? We think it is useful to keep this
picture/issue in mind.

In section~\ref{sec:sr} we briefly recall the relativity theory
$\SpecRel$ from \cite{Vienna}, in section~\ref{sec:sigth} we recall
Signalling theory $\SigTh$ from \cite{Ax} and we try to give a basic
feeling for it by sketching the proof of completeness theorem in
\cite{Ax}. Section~\ref{sec:algo} is an important part of the
present paper, it contains an algorithm for how to set up a
reference frame in Signalling theory, this is an ``operational
semantics" for setting up reference frames of \cite{Vienna}. At the
end of the section we outline how the same method could be used for
space-times other than special relativistic, e.g., for the
Schwarzshild space-time of a black hole. This algorithm is at the
heart of the interpretation elaborated in section~\ref{sec:reduc}.
Section~\ref{sec:defth} recalls the features of the more refined
definability theory that are needed for defining the interpretation
of $\SpecRel_0$ into $\SigTh$. Section~\ref{sec:defeq} rounds up the
picture between $\SpecRel$ and $\SigTh$ by interpreting $\SigTh$ in
a slightly reinforced version of $\SpecRel$ and then giving more
information about connections between various concrete theories of
special relativity. We end the paper with a Conclusion.

\section{Special Relativity}
\label{sec:sr} In this section we give a list of basic concepts and
axioms of the FOL theory $\SpecRel$ in \cite{AMNHbSL, Vienna, AMNSz,
Mad02, Sze}.

The basic notions not analyzed further in $\SpecRel$ are
``observers" having reference frames in which they represent the
world-lines of bodies (or test particles), of which signals
(light-particles, or photons) are special ones. The world-line of a
body represents its motion, it is a function that describes  the
location of the body at each instant. For representing ``time" and
``location", observers use quantities, quantities are endowed with
addition and multiplication in order to be able to express whether a
motion is ``uniform" or not. To make life simpler, we treat also
observers as special bodies. (Another, equivalent, option would be
to treat them as entities of different ``kind", or of different
``sort", than bodies and quantities.) The reference frame or
world-view of an observer $o$ gives the information which bodies $b$
are present at time $t$ at location $x,y,z$; thus $\W(o,b,t,x,y,z)$
expresses that body $b$ is present at $t$ in $\langle x,y,z\rangle$,
according to observer $o$. We treat quantities as entities of
different nature, of different kind, than bodies.

According to the above, the vocabulary of the language of $\SpecRel$
is the following: we have two sorts, bodies $\B$ and quantities
$\Q$, we have two unary relations $\Obs,\Ph$ of sort $\B$, we have
two binary functions $+,\star$ of sort $\Q$, and we have a six-place
relation $\W$ the first 2 places of which are of sort $\B$ and the
rest of sort $\Q$.

Next, we list the five axioms of $\SpecRel$. Concrete formulas and
more intuition can be found in, e.g., \cite{AMNHbSL, Vienna, AMNSz,
Mad02, Sze}.

\smallskip \noindent {\AxPh}\quad The world-lines of photons are
exactly the straight lines of slope 1, in each reference frame.

\smallskip \noindent {\AxEv}\quad All observers coordinatize the same
physical reality (i.e., the same set of events).

\smallskip \noindent {\AxSelf}\quad The ``owner" of a reference frame
sits tight (stays put) at the origin.

\smallskip \noindent {\AxFd}\quad The quantities form a Euclidean
field w.r.t.\ the operations $+,\star$, this means that $\Q,+,\star$
form an ordered field in which each positive quantity has a square
root.

\smallskip \noindent {\AxSym}\quad All observers use the same units
of measurement: if two events are simultaneous for observers $o,o'$,
then the spatial distance between them is the same according to
$o,o'$.

\bigskip
\noindent $\SpecRel_0\de\{\AxPh, \AxEv, \AxSelf, \AxFd\}$\quad and\\
\quad $\SpecRel\de\{\AxPh, \AxEv, \AxSelf, \AxFd, \AxSym\}$.
\bigskip

$\SpecRel$ may seem to be a rather weak axiom system. However, this
is not so. All the well-known theorems/predictions of the
(kinematics of) special relativity can be proved even from
$\SpecRel_0$. Below is a sample of theorems that can be proved from
$\SpecRel_0$ (for proofs, further theorems provable from $\SpecRel$,
and for extensions see the references given earlier as well as
\cite{MadSzek13,Sze10}):

\begin{description}
\item[$\bullet$]
Each observer moves uniformly and slower than light in any other
observer's world-view (i.e., the world-line of an observer is a
straight line with slope less than 1).
\item[]
Assume that $o,o'$ are moving relative to each other.
\item[$\bullet$]
Events that are separated in $o$'s world-view in a direction of
$o'$'s motion and simultaneous according to $o$, are not
simultaneous according to $o$'.
\item[$\bullet$]
Events that are simultaneous according to both $o$ and $o'$ are
exactly the ones that are separated orthogonally to the direction of
motion of $o'$.
\item[$\bullet$]
Assume that $o$ and $o'$ use the same units of measurement, i.e.,
the spatial distance between events that are simultaneous to both of
them is the same according to them. Then a-synchronicity,
time-dilation and length-contraction between $o$ and $o'$ are
exactly according to the known formula of special relativity, see
e.g., \cite[p.633]{AMNHbSL}.
\item[$\bullet$]
The world-view transformations between observers in $\SpecRel_0$ are
exactly the bijections that preserve Minkowski-equidistance; these
bijections are the so-called Poincar\'e-transformations composed
with dilations and field-automorphisms.
\item[$\bullet$]
The world-view transformations between observers in $\SpecRel$ are
exactly the bijections that preserve Minkowski-distance; these
bijections are the so-called Poincar\'e-transformations.
\end{description}
\bigskip

In $\SpecRel$, a reference frame is a basic (or primitive) notion,
just an ``out-of-the-blue" assigning space-and-time coordinates to
events, which all together have to satisfy some regularities (our
axioms). The theory does not address the question of how an observer
sets up his reference frame. As already outlined in the
introduction, according to some authors, a physical theory (a theory
about our physical reality), should say something about the meaning
(in the ``real" physical world) of the basic concepts, if not
otherwise, then in natural language one could amend the theory with
a set of so-called \textit{operational rules} about how the basic
concepts (the reference frames in our case) are set up
(experimentally). Here usually meter-rods and wrist-watches, or
standard clocks, are used, see e.g., Taylor and Wheeler \cite[Fig.s
9,135]{Wheeler}, L.\ E.\ Szabo \cite{Szabo}. In
section~\ref{sec:algo} we give a more ambitious algorithm for
setting up coordinate systems.

\section{James Ax's Signalling theory} \label{sec:sigth}

The intention of Ax's theory is to give an axiom system for special
relativity so that its basic symbols and axioms are designed to be
observational. The players of this theory are experimenters that can
``communicate with each other" by sending signals to each other.
Together, as a team the experimenters can ``map" (or explore)
space-time, without having rigid meter rods or clocks. A definition
of an introduced (or defined) term in this first-order logic theory
can be viewed as an experiment designed to establish whether the
defined term holds or not. The basic terms of space and time are
defined this way. (Indeed, in this theory one can define ``rigid
meter rods" and ``clocks" from signalling experiments.) The results
of the experiments we make can be built into axioms then (which are
designed to be observational-oriented), and they can tell us in what
kind of space-time we live in. Euclidean? Special relativistic?
Hyperbolic space with relativistic time? Newtonian? General
relativistic?  Etc. All this amounts to an implementation of
Leibnitzian relational notion of space and time. We return to this
subject in more detail in the next section.

We begin to describe Ax's theory which we call Signalling theory
$\SigTh$. In the vocabulary of $\SigTh$ we have two sorts, $\Par$
for ``particles" (or experimenters, or agents) and $\Sig$ for
``signals" (or light-signals); and we have two binary relations
$\T,\R$ between particles and signals. The intended meanings of
$a\T\s$ and $a\R\s$ are ``$a$ transmits (or emits, or sends out)
$\s$", and ``$a$ receives (or absorbs) $\s$", respectively. Ax
\cite{Ax} uses an impersonal terminology of particle physics,
particles emit and absorb signals. We are more attracted to a
terminology of communication between active experimenters. These
experimenters (players) of $\SigTh$ are somewhat analogous to the
observers of $\SpecRel$. In this paper when talking about Ax's
Signalling theory, we will use the terms experimenter and particle
interchangeably.

The ``standard" (or intended) model we have in mind is the
following: Let us fix a Euclidean field $\F$. Then $\Par$ is the set
of all straight lines in $\F^4$ with slope less than 1, and $\Sig$
is the set of all directed finite segments (including the segments
of length zero) of straight lines with slope 1. A particle $a$
transmits a signal $\s$ iff the beginning point of $\s$ lies on $a$,
and $a$ receiving the signal $\s$ means that the endpoint of $\s$
lies on $a$. Let us denote this structure by $\Mm(\F)$.

The main result of \cite{Ax} is a finite set $\Sigma$ of axioms, our
$\SigTh$, which characterizes the class of standard models, i.e.,
the models of $\Sigma$ are exactly the standard models $\Mm(\F)$
over some Euclidean field $\F$ (Thm.1 in \cite{Ax}). $\SigTh$
consists of three groups of axioms, altogether it has 23 elements.
Instead of listing these 23 axioms, in this paper we will use Ax's
completeness theorem, since that implies that a formula $\psi$ is
provable from $\SigTh$ iff $\psi$ is true in all the standard
models.

The question immediately arises: what does this theory $\SigTh$ have
to do with special relativity? Do not we lose much expressive power
by using such meager resources? Where do a-synchronicity,
time-dilation, length-contraction come into the picture in $\SigTh$?
Some answers are in the proof of Thm.1 which we briefly outline
below. We will give more explicit answers to these questions in the
coming sections~\ref{sec:algo}, \ref{sec:reduc}. In particular, we
will show that everything we can say in the language of $\SpecRel$
can be said in the Spartan language of $\SigTh$, too. One of the
ideas for proving this can be traced back to Hilbert, as will be
noted below Figure~\ref{times-fig}.

To give a feeling for $\SigTh$ and the expressive power of its
language, we briefly outline the proof of Ax's completeness theorem.
Let's begin by making a little elbow-room for working. We will need
to express things such as ``two signals are received by an
experimenter at the same time", and ``signal $\s$ was received by an
experimenter just when he transmitted signal $\g$". Since we have no
notion of time in our language, we have to express these notions
just by using the basic concepts of transmitting and receiving
signals. Here comes how we can do this. The (open) formula
``$\phi\de\forall a\enskip a\T\s\to a\T\g$" is true in a standard
model just when the beginning points of the segments $\s,\g$
coincide, we say that $\phi$ expresses this fact.%
\footnote{In the formulas, the scope of a quantifier is till the end
of the formula if not indicated otherwise. Lower case Roman and
Greek letters denote variables of sorts $\Par$ and $\Sig$,
respectively. In place of conjunction $\land$ we will simply write a
comma.} Similarly, ``$\forall a\enskip a\R\s\to a\T\g$" expresses
that the endpoint of $\s$ coincides with the beginning point of
$\g$, etc. Now can we express that two particles/experimenters meet?
Well, they meet if there is a signal that both of them transmit.
From now on we will use similar statements without translating them
to the language of $\SigTh$.

To begin outlining the idea of Ax's completeness proof, let $\Mm$ be
any model of $\SigTh$, and let $\ea$ be any experimenter in this
model. We will construct an isomorphism between $\Mm$ and a standard
model $\Mm(\F)$ which takes $\ea$ to the time-axis in $\Mm(\F)$.
From now on, in this section $\ea$ denotes this fixed experimenter.

We define the set $\Space$ of ``places" or ``locations" for
experimenter $\ea$ to consist of those particles which are
motionless w.r.t.\ $\ea$. For expressing that two particles are
motionless w.r.t.\ each other, any formula expressing this in the
standard models will do. Ax uses the following formula: $\ea'$ is
motionless w.r.t.\ $\ea$ exactly when $\ea$ and $\ea'$ do not meet
(if they are not equal) and there are two other particles $d,c$
which meet them and each other in 5 distinct events. Ax then
expresses the betweenness relation $\Bw$ for such places as well as
the equidistance relation $\Ed$ with suitable formulas. Having all
this, the first group of axioms in $\SigTh$ states Tarski's axioms
for axiomatizing Euclidean geometry over the Euclidean fields (see
\cite{Tarski}). From Tarski's theorem then Ax gets a Euclidean field
$\F$ and an isomorphism between $\langle\F^3,\Bw,\Ed\rangle$ and
$\langle\Space,\Bw,\Ed\rangle$.

Having $\Space$ for our experimenter $\ea$, what is ``time" for it?
What are the things that we mark with time? The events. And what are
the events? In the present vocabulary we take them to be ``particle
$b$ emits/receives a signal $\s$", more precisely we take the
equivalence classes of them described when we made the elbow-room
for this proof (e.g., particle $b$ may send out signal $\s$ in the
same event when it sends out another signal $\g$ or when it receives
$\g$). Then our experimenter $\ea$'s time will be the events that
happened to $\ea$. For simplicity, we will represent events with
special signals, as explained below.

In the standard models, there are special signals that are received
by everyone who transmitted them, we call these signals {\it
events}, we will denote them by variants of $\e$:

\begin{description}
\item{} $\Ev(\e)\deiff \forall a\enskip a\T\e\to a\R\e$.
\end{description}

\noindent In the standard models, events are the light-like segments
of zero length, so they correspond to elements of $\F^4$. (These
zero-length signals may look counter-intuitive to some readers. It
is just handy and not important that we use or have these at hand,
everything works with a slight modification if we omit these short
signals from the standard models.) We say that event $\e$ happened
to experimenter $\ea$, or in other words, experimenter $\ea$ {\it
participated in event} $\e$, if $\ea$ transmitted (and then also
received) $\e$. The events that happened to $\ea$ will constitute
$\ea$'s {\it world-line}.

We can then express simultaneity of events by using that the speeds
of light-signals are the same (see Figure~\ref{sim-fig}). Ax then
states an axiom to the effect that signals make a one-to-one
correspondence between $\ea$'s world-line and the simultaneous
events on any given line in $\Space$. This makes $\ea$'s world-line
isomorphic to $\F$, we take this to be the time-axis.
From now on it is more or less straightforward what we have to
include to $\SigTh$ in order to make $\Mm$ isomorphic to $\Mm(\F)$.
E.g., we can state that for any event $\e$ there is a simultaneous
event $\e'$ on the world-line of $\ea$, and there is a particle
$\ea'$ that participates in $\e'$ and is motionless w.r.t.\ $\ea$.
This concludes the proof-idea.

\section{An algorithm for setting up coordinate systems} \label{sec:algo}

The purpose of this section is twofold. Firstly, in
section~\ref{sec:reduc} we want to give an interpretation of
$\SpecRel_0$ into $\SigTh$, and for this we need concrete formulas
representing the proof-idea given in the previous section. For
example, Ax used Tarski's theorem for getting the Euclidean field
$\F$, but we will need to exhibit concrete formulas defining this
$\F$. Secondly, we want to make the previous proof-idea into an
algorithm for setting up coordinate systems (i.e., reference frames)
with the use of just light-signals and freely moving particles. This
could also be viewed as providing {\it operational semantics} to the
basic notion of a coordinate system of $\SpecRel$. What we give in
this section will not be an algorithm in the strict sense, it will
be more like a recipe for how to design experiments/measurements for
assigning coordinates to events. These experiments will also be
suitable for finding out/confirming that we live in a special
relativistic space-time (if we do). For this reason, we will try to
make the formulas ``executable" when possible. There will be plenty
of room for improving on this aspect, the reader is invited to
design more practical experiments.

Assume that we are given a model $\Mm$ of $\SigTh$, and $\ea$ is an
experimenter in this model. Just as in the previous section, this
experimenter $\ea$ is fixed throughout this section. We are going to
give $\ea$ a recipe for defining a field $\F$ of quantities and for
assigning four quantities to each event. Such an assignment is
called a {\it coordinate system} (or reference frame). These
coordinate systems will satisfy the axioms of $\SpecRel_0$.

A {\it location} for $\ea$ was defined as a particle that is
motionless w.r.t.\ $\ea$. In the previous section we recalled a
formula, from \cite{Ax}, expressing whether $\ea'$ is motionless
w.r.t.\ $\ea$ (in symbols, $\ea'\| \ea$). However, the algorithm
suggested by that formula is not very convenient since it involves
deciding whether $\ea$ meets $\ea'$ or not, and for this $\ea$ has
to know all the events that happened and will happen to him. This is
not very practical as an experiment, since $\ea$ may need to ``wait"
for an infinity of time before he could know the result. Using the
Affine Desargues Property (ADP for short, see, e.g.,
\cite[p.20]{Gold}) one can design a more realistic experiment which
decides $\ea'\| \ea$ ``in a finite time", we are going to describe
it now. We note that in the standard models $\Mm(\F)$ the ADP is
true, because it is true in the affine space $\F^4$, for any field
$\F$.

\begin{figure}[h!]
\begin{center}
\psfrag{a}[b][b]{$\ea$} \psfrag{a'}[b][b]{$\ea'$}
\psfrag{e}[b][b]{$\e$}\psfrag{e1}[b][b]{}\psfrag{e2}[b][b]{}
\psfrag{e'}[b][b]{$\e'$} \psfrag{?}[m][b]{\large ?}
\psfrag*{b1}[b][b]{$b_1$}\psfrag*{b2}[b][b]{$b_2$}\psfrag*{b3}[b][b]{$b_3$}
\includegraphics[scale=.3]{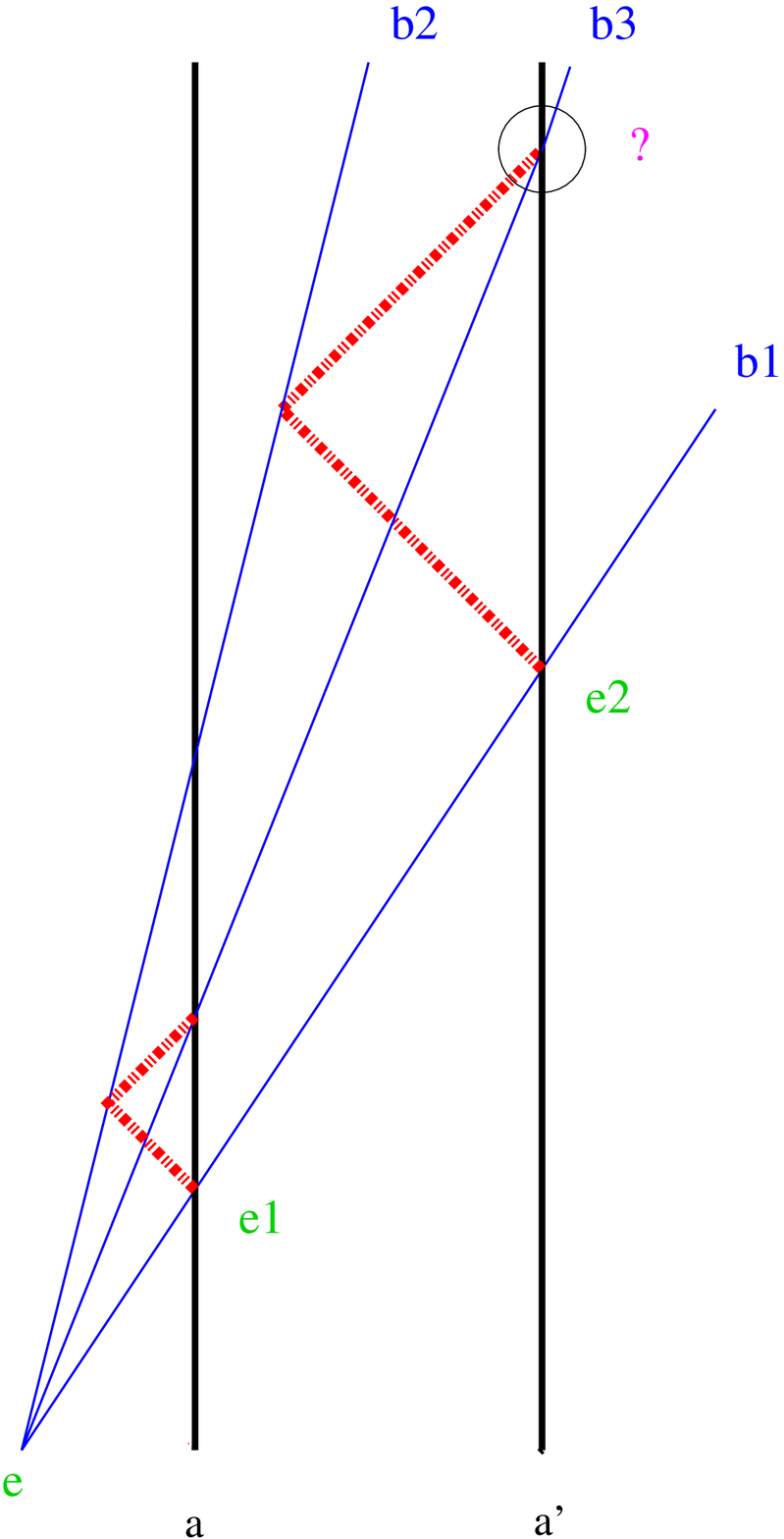}\hskip 1.5cm
\psfrag{a}[b][b]{$\ea'$}
\psfrag{e1}[b][b]{$\e_1$}\psfrag{e2}[b][b]{$\e_2$}
\includegraphics[scale=.25]{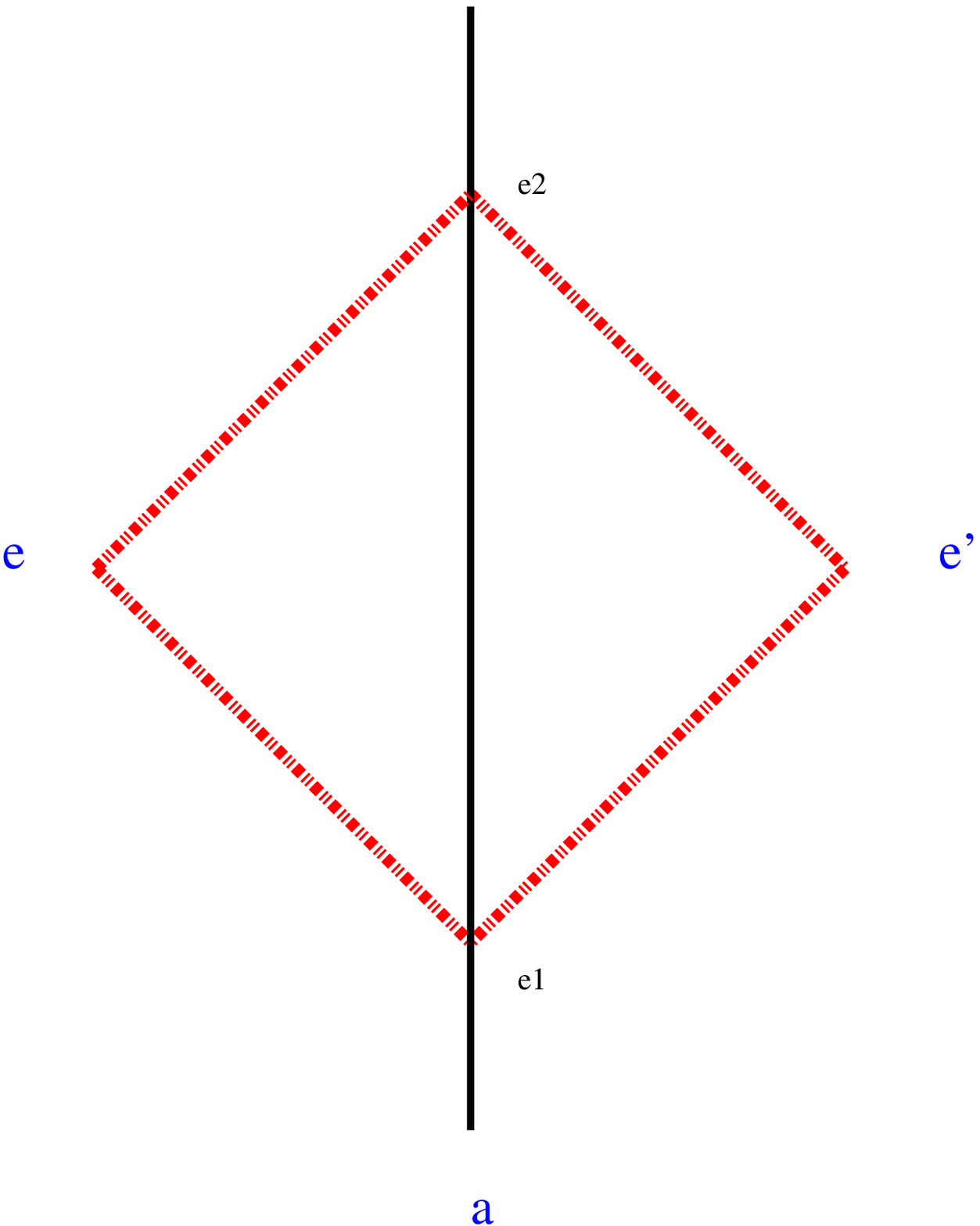}\hskip 2.1cm
\psfrag{a}[b][b]{$\ea$}
\psfrag{e1}[b][b]{$\e_1$}\psfrag{e2}[b][b]{$\e_2$}
\psfrag{e3}[b][b]{$\e_3$}\psfrag{e4}[b][b]{$\e_4$}
\includegraphics[scale=.3]{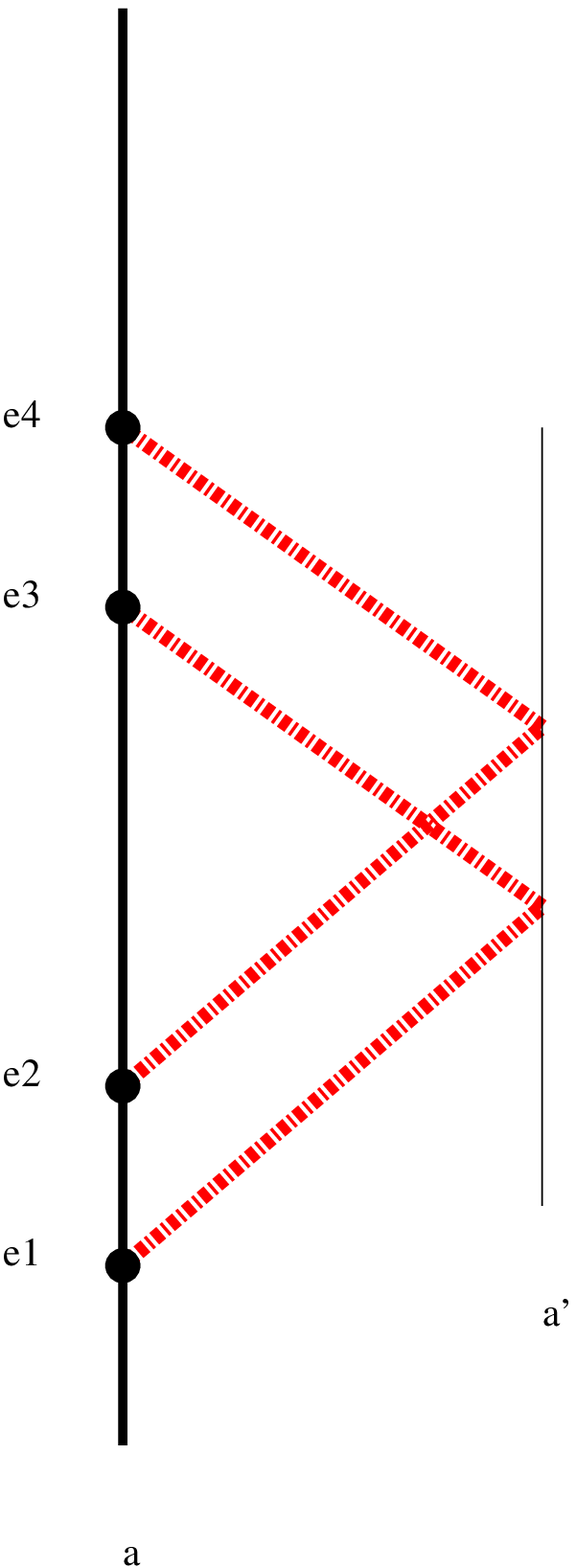}
\end{center}
\caption{\label{des-fig}\label{sim-fig} On the left: Experiment for
checking whether $\ea'$ is motionless w.r.t.\ $\ea$. In the middle:
Experiment to make sure that $\e,\e'$ are simultaneous w.r.t.\
$\ea$. On the right: Time-equidistance of events $\e_1,\dots,\e_4$.}
\end{figure}

For a while, it will be easier to think in 4-dimensional space-time
than tracing motion in 3-dimensional space. Geometrically, $\ea'$ is
motionless w.r.t.\ $\ea$ iff the world-line of $\ea'$ is parallel to
that of $\ea$. The conclusion of the ADP is that two lines are
parallel, but in the hypothesis part parallelism of two other sets
of lines are used. We are lucky: we have light-signals and their
speeds are the same in both directions, thus we can use parallelism
of world-lines of two sets of light-signals in the hypothesis part
of the ADP. The experiment is depicted in geometrical form in the
left-hand part of Figure~\ref{des-fig}. Here is the
``non-geometrical" description of the experiment: Assume $\ea$ wants
to decide whether $\ea'$ is motionless w.r.t.\ him or not. He  asks
a brother (another experimenter) to throw towards him three ``test"
particles (``balls") $b_1,b_2,b_3$ at once (in one event $\e$),
$b_1$ faster than $b_3$ and $b_3$ faster than $b_2$ in such a way
that when $b_1$ meets $\ea$, the latter sends out a signal towards
$b_2$ that $b_2$ reflects back and the reflected signal reaches
$\ea$ just when $b_3$ reaches $\ea$. (The brother and $\ea$ have to
experiment a little while till finding the right velocities for such
three particles.) After checking that $b_1,b_2,b_3$ have the desired
property, $\ea$ asks $\ea'$ to do the same: when $b_1$ reaches $a'$,
he should send a signal towards $b_2$ that reflects this signal back
towards $\ea'$. If the reflected signal reaches $\ea'$ just when
$b_3$ reaches $\ea'$, then $\ea'$ is motionless w.r.t.\ $\ea$;
otherwise $\ea'$ is not motionless w.r.t.\ $\ea$. It is best to
imagine this experiment to take place in outer space, far from heavy
heavenly objects so that gravity and
friction do not bend the world-lines of the ``balls".%
\footnote{Or, if we are content with more approximate measurements,
we can imagine all this happening on a big lake covered with smooth
ice (but then we have to take space to be 2-dimensional).}
 From now on, we will use ``locations" and ``places" as being particles/experimenters
  motionless w.r.t.\ our fixed experimenter $\ea$.

Two events $\e,\e'$ are defined to be {\it simultaneous} w.r.t.\
$\ea$ iff there is a place $\ea'$ such that from $\ea'$ two signals
can be sent at the same event towards the locations of $\e$ and
$\e'$ respectively such that if these signals are sent back from
$\e$ and $\e'$ right away, they will arrive back to $\ea'$ at the
same event, see middle of Figure~\ref{sim-fig}. Formally:
$\e\simu_{\ea}\e'\deiff \\
\exists \ea'\| \ea, \s_1,\dots,\s_4, \e_1,\e_2\enskip \Ev(\e_1),
\Ev(\e_2), (\e_1,\s_1,\e), (\e_1,\s_2,\e'), (\e,\s_3,\e_2),
(\e',\s_4,\e_2), \ea'\Tran\e_1, \ea'\Tran\e_2$,\ where $(\e,\s,\g)$
means that $\e,\g$ are the events of sending and receiving $\s$,
respectively, formally: $(\e,\s,\g)\deiff \\
(\Beg(\s,\e), \End(\s,\g))$ where $\Beg(\s,\e)$ expresses that
$\s,\e$ are sent out at the same event, formally: $\Beg(\s,\e)\deiff
\forall b\enskip b\T\s\to b\T\e$ and a similar definition for
$\End$. (Note that if we want a more experiment-friendly formula for
$\Beg$, then we can use the following:
$\Beg(\s,\e)\Leftrightarrow(\exists b,c\enskip b\ne c, b\T\s,
b\T\e,\\ c\T\s, c\T\e)$.) We even can provide instructions for where
to look for such a place $\ea'$: it can be chosen to be the midpoint
of the line-segment connecting the locations of $\e$ and $\e'$. (We
can use this experiment for setting two clocks at the places of
$\e,\e'$ which ``tick simultaneously".)

We get an ordering on all the events from the fact that we send a
signal earlier than receiving it, namely $\e$ is {\it earlier} than
$\e'$ iff we can send a signal at $\e$ to an event from where it
bounces back to $\e'$ ($\e\prec\e'\deiff [\exists
\e'',\s_1,\s_2\enskip (\e,\s_1,\e''), (\e'',\s_2,\e')]$). For
example, $\e_1$ is earlier than $\e_2$ in the middle part of
Figure~\ref{sim-fig}. We note that, while two events being
simultaneous or not depends on which experimenter makes the
experiment deciding simultaneity, one event being earlier than
another does not depend on any experimenter.

Let's see, what structure the set of events happening to $\ea$ has.
Let $\Time_{\ea}$ denote the world-line of $\ea$ and let
$\e_1,\dots,\e_4\in\Time_{\ea}$. Besides the ordering, we also have
{\it time-equidistance} of events, since the speed of all signals is
the same: the time elapsed between $\e_1$ and $\e_2$ is the same as
that between $\e_3$ and $\e_4$, iff there is a place $\ea'$ to which
we can send signals from $\e_1,\e_2$ resp., these bounce from $\ea'$
and arrive back to $\ea$ at $\e_3,\e_4$ respectively. See the
right-hand part of Figure~\ref{des-fig}. More precisely, this is the
definition  when $\e_1\prec\e_3$.  When $\e_3\prec\e_1$, we get the
definition by interchanging the pairs $\e_1,\e_2$ and $\e_3,\e_4$.
(Formally,
$\Ted_{\ea}(\e_1,\e_2,\e_3,\e_4),\mbox{$\e_1\prec\e_3$}\deiff
[\exists \ea'\| \ea, \s_1,\dots,\s_4,\e,\e'\enskip \Ev(\e),
\Ev(\e'),\\ (\e_1,\s_1,\e), (\e,\s_3,\e_3), (\e_2,\s_2,\e'),
(\e',\s_4,\e_4), \ea'\T\e, \ea'\T\e']$,\ \  and\ \
$\Ted_{\ea}(\e_1,\e_2,\e_3,\e_4),\e_3\prec\e_1\deiff\\
\Ted_{\ea}(\e_3,\e_4,\e_1,\e_2),\e_3\prec\e_1$.)  Note that
$\Ted_{\ea}(\e_1,\dots,\e_4)$ implies that $\e_1$ happens earlier
than $\e_2$ iff $\e_3$ happens earlier than $\e_4$. By using
time-equidistance, we can define {\it addition} by selecting a
``{\it zero}" time $\oo\in\Time_{\ea}$ as parameter, namely
$\st=\st_1+\st_2\deiff +\!(\st,\st_1,\st_2,a,\oo)\deiff
\Ted_{\ea}(\oo,\st_1,\st_2,\st)$. Now that we have addition, we do
not stop before having {\it multiplication}. For this we have to
choose a {\it unit} time $\uu\in\Time_{\ea}$, distinct from $\oo$
and happening later than $\oo$, as another parameter. For defining
multiplication, we will need the {\it collinearity} relation on
locations, we will get this by noticing that the space-trajectories
of signals are (3-dimensional) straight lines in the standard
models: $\Col(a_1,a_2,a_3)$ iff exist signals $\s_1,\s_2,\s_3$ and
events $\e_1,\e_2,\e_3$ such that $(\e_i,\s_1,\e_j),
(\e_j,\s_2,\e_k), (\e_i,\s_3,\e_k)$ and $a_1\T\e_1$, $a_2\T\e_2,
a_3\T\e_3$, for some permutation $i,j,k$ of $1,2,3$.

\begin{figure}[h!]
\begin{center}
\psfrag{t1*t2}[b][l]{$\st_1\star \st_2$} \psfrag{t2}[b][l]{$\st_2$}
\psfrag{t2'}[t][t]{$\st_2'$} \psfrag{t1}[b][l]{$\st_1$}
\psfrag{u}[b][l]{$\uu$}\psfrag{u'}[t][t]{$\uu'$}
\psfrag{o}[b][b]{$\oo$} \psfrag{b}[b][b]{$b$}
\psfrag{az}[b][b]{$a_z$} \psfrag{ax}[b][b]{$a_x$}
\psfrag{ay}[b][b]{$a_y$} \psfrag{ex}[b][b]{$p_x$}
\psfrag{ey}[b][b]{$p_y$} \psfrag{ez}[b][b]{$p_z$}
\includegraphics[scale=.35]{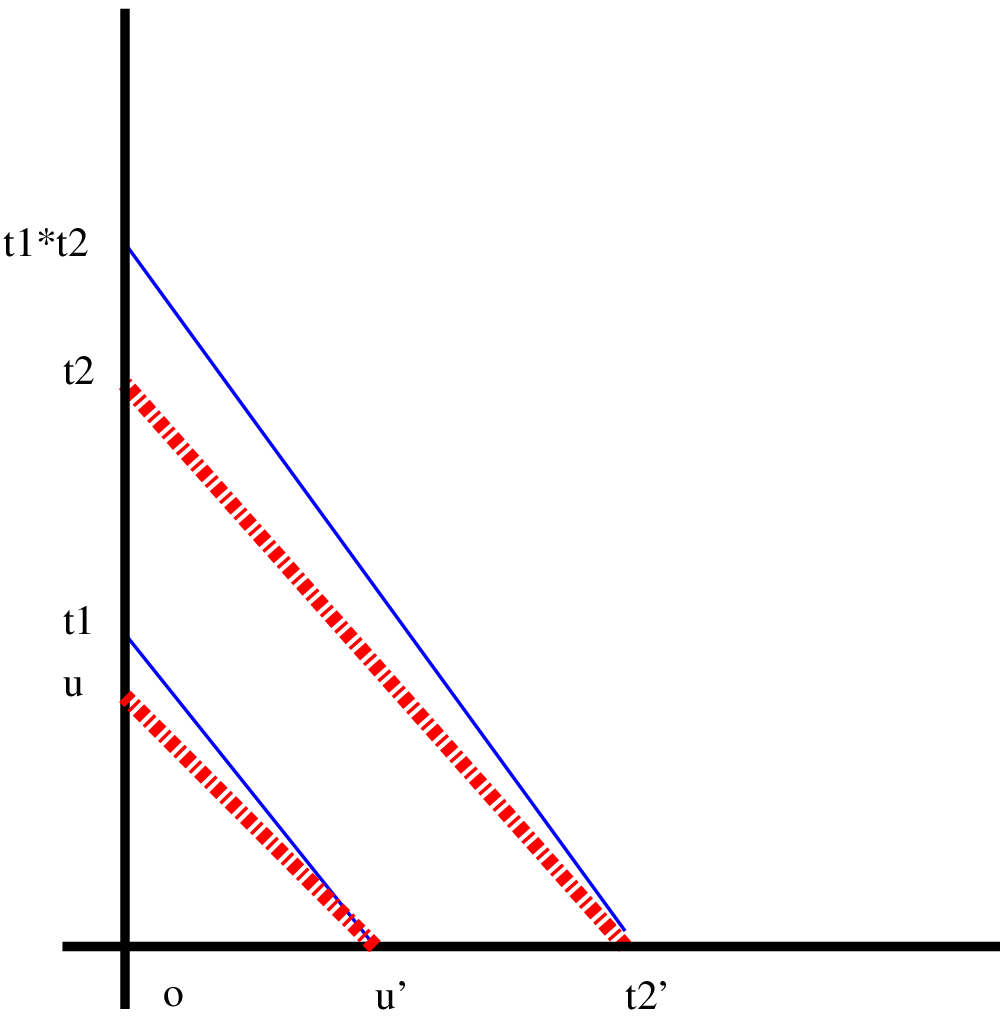}\hskip 1.5cm
\psfrag{ee}[t][t]{$\equiv_{\ea}$} \psfrag{dab}[r][r]{$\dd(\ea,b)$}
\psfrag{a}[b][b]{$\ea$}
\includegraphics[scale=.3]{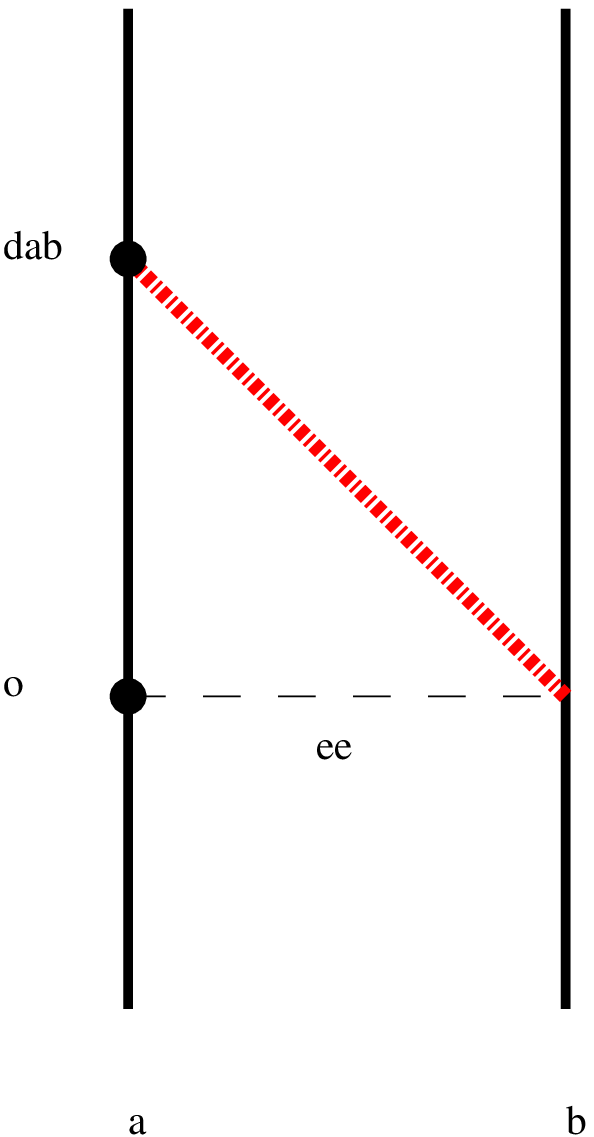}\hskip 2.5cm
\psfrag{az}[b][b]{$a_z$} \psfrag{ax}[b][b]{$a_x$}
\psfrag{ay}[b][b]{$a_y$} \psfrag{gx}[b][b]{$\g_x$}
\psfrag{o}[b][b]{$\ea$} \psfrag{b}[b][b]{$b$}
\psfrag{ez}[b][b]{$p_z$} \psfrag{ex}[b][b]{$p_x$}
\psfrag{ey}[b][b]{$p_y$}
\includegraphics[scale=.2]{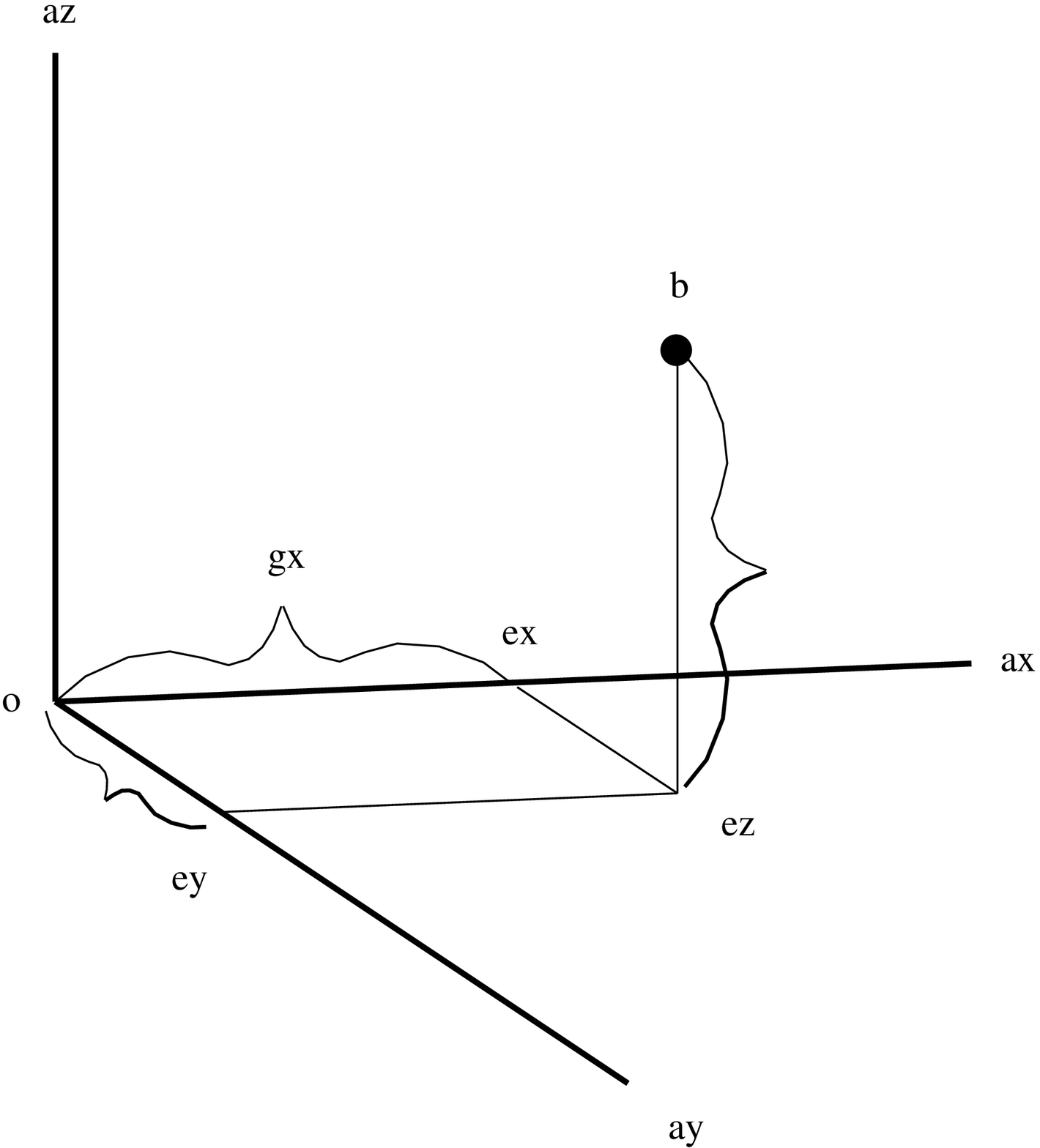}
\end{center}
\caption{\label{times-fig}\label{coord-fig} On the left: Experiment
for computing $\st_1\star \st_2$. In the middle: distance between
locations $\ea,b$. On the right: Spatial coordinates of location
$b$. In the picture, $\g_x=\dd(\ea,p_x)$. In this part points
represent (3-dimensional) locations, while in the previous pictures
points represent (4-dimensional) events.}
\end{figure}

We define $\st_1\star \st_2$ for the case when $\st_1$ happened
later than $\uu$, and $\st_2$ happened later than $\oo$. See the
left-hand part of Figure~\ref{times-fig}. (The other cases are
similar, we leave them out.) Here is how we find out whether $\st$
is $\st_1\star \st_2$: we find two places $b_1$ and $b_2$ collinear
with $\ea$ and we find a particle $p$ such that if $b_1$ and $b_2$
send towards $\ea$, simultaneously, at time zero, a light-signal and
$p$, and another light-signal and another particle $q$ ``with the
same speed" as $p$, then these four arrive (to $\ea$) at times $\uu,
\st_1, \st_2, \st$, respectively. Formally: $\mbox{$\st=\st_1\star
\st_2$}\deiff \mbox{$\star(\st,\st_1,\st_2,\ea,\oo,\uu)$}\deiff
\exists b_1\| \ea, b_2\| \ea, \uu',\st_2', \s_1,\s_2,\\ p,q\enskip
[\Ev(\uu'), \Ev(\st_2'), \Col(\ea,b_1,b_2),
 \oo\simu_{\ea} \uu', \oo\simu_{\ea} \st_2', b_1\T \uu', b_2\T \st_2', (\uu',\s_1,\uu),
 (\st_2',\s_2,\st_2), p\T \uu', p\T \st_1, q\T \st_2',\\ q\T \st, p\| q]$.

The reader will have noticed that the above definitions of addition
and multiplication on $\Time_{\ea}$ are a special case of Hilbert's
coordinatization procedure, see e.g., \cite[pp.23-28]{Gold} or
\cite[pp.296-308]{Mad02}.

By the above, we have a structure $\F(\ea,\oo,\uu)=\langle
\Time_{\ea}, +, \star\rangle$ which is isomorphic to our field $\F$
in the intended models $\Mm(\F)$. We define the above structure to
be the field of {\it quantities} of our fixed experimenter $\ea$. We
define the {\it time-coordinate} of an arbitrary event $\e$ as an
element of this field, namely the unique event on $\ea$'s world-line
which is simultaneous with it (simultaneous according to $\ea$).
Next, we define three coordinates, three elements of this field, for
each location $b$. From now on, let $\Space_{\ea}$ denote the set of
locations for $\ea$.

We begin by defining a geometric structure on $\Space_{\ea}$, namely
we will define distance of locations, parallelism and orthogonality
of (3-dimensional) spatial lines.

We define the {\it distance} of any two locations. Let
$b\in\Space_{\ea}$ be arbitrary. We define the distance of $b$ from
our fixed $\ea$ as the event when a signal sent from $b$ at time
zero arrives to $\ea$, see the middle part of
Figure~\ref{times-fig}. This definition corresponds to a convention
that we measure spatial distances in light-years (if we measure time
in years). Having this, we get the distance between any two
locations $b_1,b_2$ by measuring the distance between their parallel
translated versions so that $b_1$ gets to $\ea$, i.e.,
$\dd(\ea,b)=\e\deiff \exists\e',\s\,[ \Ev(\e'), \e'\simu_{\ea} \oo,
b\T\e', (\e',\s,\e), \e\in\Time_{\ea}]$, and
$\dd(b_1,b_2)=\e\deiff[\e=\dd(\ea,b), \pa(b_1,b_2,\ea,b),
\pa(b_1,\ea,b_2,b)]$ where $\pa(b_1,b_2,b_3,b_4)$ means that the
spatial lines defined by $b_1,b_2$ and $b_3,b_4$ are {\it parallel},
we easily can express this by using the collinearity relation $\Col$
between locations as defined earlier in this section.

We also need the orthogonality relation which is definable from the
equidistance of pairs of locations. We define orthogonality of two
intersecting lines only. We call the lines going through $a,b$ and
$a,c$ orthogonal, if $a\ne b, a\ne c$ and there is a $b'\ne b$ on
the spatial line going through $a,b$ such that the distances between
$a,b'$ and $a,b$ equal, and also those between $c,b'$ and $c,b$
equal ($\Ort(a,b,a,c)\deiff \exists b'[\Col(b',a,b),
\dd(a,b')=\dd(a,b), \dd(c,b')=\dd(c,b)]$). By now we defined a
structure $\langle\Space_a,\Col,\pa,\Ort\rangle$ and we defined
distance $\dd:\Space_{\ea}^2\longrightarrow\Time_{\ea}$.

Setting up a coordinate system needs three more parameters, the
three space-axes. Let $a_x, a_y, a_z\in\Space_{\ea}$ be such that
$\ea,a_x$, $\ea,a_y$ and $\ea,a_z$ are pairwise orthogonal. We have
everything for defining the usual spatial coordinates of the place
$b$. See the right-hand part of Figure~\ref{times-fig}. The spatial
coordinates of a location $b$ are defined the usual way by
``projecting" $b$ to the three coordinate axes, along lines parallel
with some of the axes, and measuring the distance of the projected
points from the origin (our experimenter $\ea$ in our case). See the
formula $\cord$ below.

We can now round up the definition of the coordinate system our
experimenter $\ea$ is setting up. We already defined the
time-coordinate of an event $\e$, and we define the
space-coordinates of $\e$ to be the spatial coordinates just defined
for the ``location of $\e$", the latter being the unique particle
participating in $\e$ and motionless w.r.t.\ our experimenter $\ea$.
The formula $\cord(\e,\st,\g_x, \g_y, \g_z,
\ea,\oo,\uu,a_x,a_y,a_z)$ defined below expresses that the
coordinates of the event $\e$ are $\st,\g_x, \g_y, \g_z$ in the
coordinate system specified by $\ea,\oo,\uu,a_x,a_y,a_z$.

\begin{description}
\item[]
$\cord(\e,\st,\g_x, \g_y, \g_z, \ea,\oo,\uu,a_x,a_y,a_z)\deiff
\e\equiv_{\ea}\st, \exists b,p_x,p_y,p_z\in\Space_{\ea} [b\T\e,
\\\pa(b,p_z,\ea,a_z), \pa(p_z,p_x,\ea,a_y), \pa(p_z,p_y,\ea,a_x),
\Col(\ea,p_x,a_x),\Col(\ea,p_y,a_y),\\ \dd(p_z,p_x)=\g_x,
\dd(p_z,p_y)=\g_y, \dd(b,p_z)=\g_z.]$.
\end{description}

\noindent Since it can be proved that the associated coordinates are
unique, we will also use the functional form

\begin{description}
\item[]
$\cord(\e,\ea,\oo,\uu,a_x,a_y,a_z)=(\st,\g_x, \g_y, \g_z)\deiff
\cord(\e,\st,\g_x, \g_y, \g_z, \ea,\oo,\uu,a_x,a_y,a_z)$.
\end{description}

\begin{figure}[h!]
\begin{center}
\psfrag{a}[b][b]{$\ea$} \psfrag{a'}[b][b]{$\ea'$}
\psfrag{a''}[b][b]{} \psfrag{e1}[b][b]{$\e_1$}
\psfrag{e2}[b][b]{$\e_2$} \psfrag{e}[b][b]{$\e$}
\psfrag{o}[b][b]{$\oo$}
\psfrag{mu}[r][r]{$\mu_{\ea,\so}(\oo,\e)=\xi$}
\psfrag{s1}[m][m]{$\equiv_{\ea}, \equiv_{\ea'}$}
\psfrag{e'}[b][l]{$\e'$} \psfrag{e''}[b][b]{$\e''$}
\includegraphics[scale=.3]{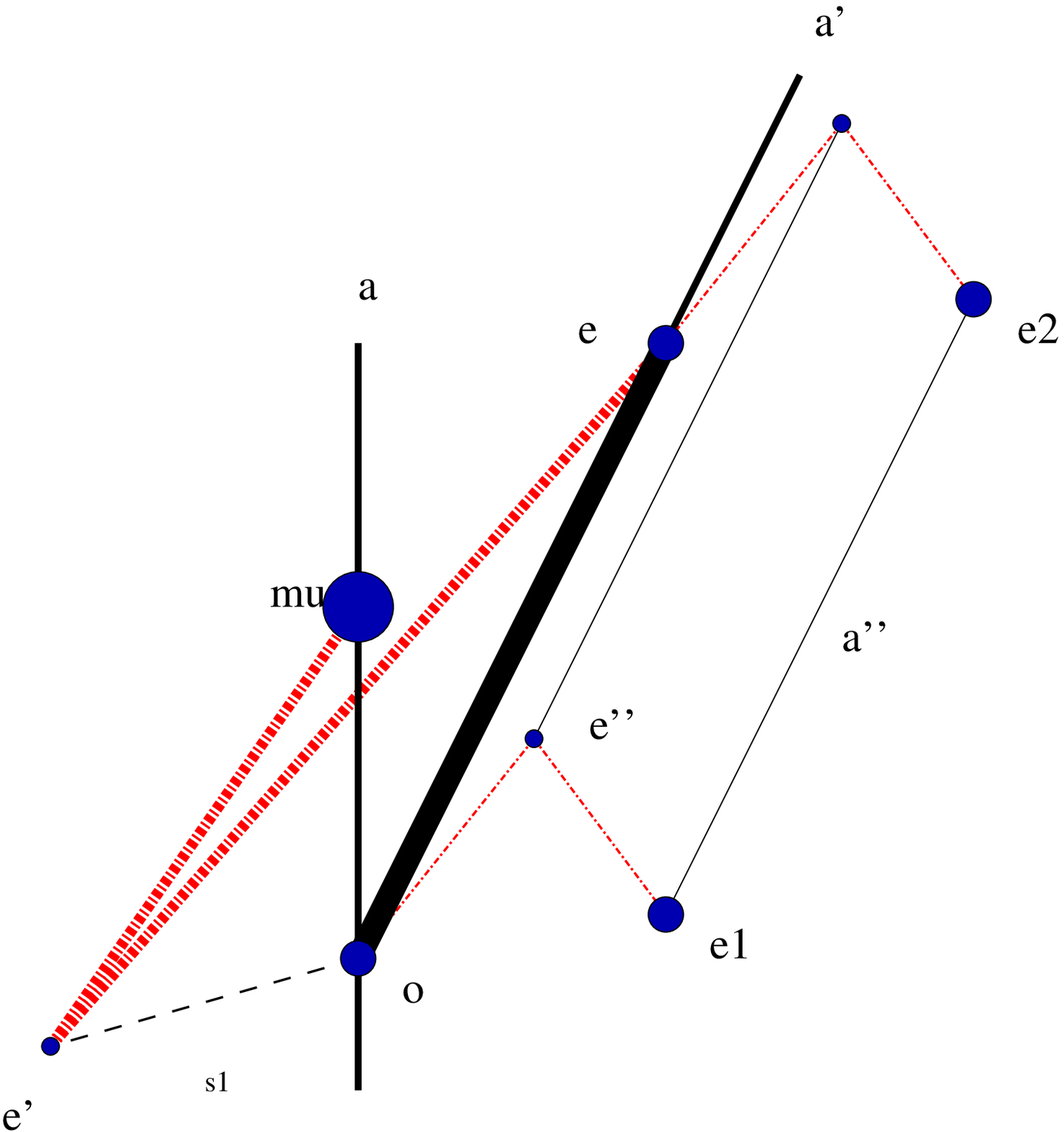}
\end{center}
\caption{\label{med-fig} Relativistic distance
$\x=\mu_{\ea,\so}(\oo,\e)=\mu_{\ea,\so}(\e_1,\e_2)$ between events
$\e_1,\e_2$.}
\end{figure}

By the above, we have defined coordinate systems to each particle
$\ea\in\Par$. Such a coordinate system is defined by six parameters:
$\ea,\oo,\uu,a_x,a_y,a_z$. Before going on, we show that the
relativistic (or, in other words, Minkowski-) distance between
events can be defined in these coordinate systems. We call two
events $\e,\e'$ {\it time-like separated} iff there is a particle
participating in both. For simplicity, we will define relativistic
distance between time-like separated events only. See
Figure~\ref{med-fig}.

The relativistic distance we are going to define will depend on
experimenter $\ea$ and on the chosen zero $\oo$ of its coordinate
system. Let first $\e\succ\oo$ be any event time-like separated from
$\oo$. Then $\mu_{\ea, \so}(\oo,\e)=\xi$ iff there is an event $\e'$
which is simultaneous with $\oo$ both according to $\ea$ and
according to the unique observer participating in $\oo,\e$, and
there are signals from $\e'$ to $\e$ and from $\e'$ to $\xi$,
respectively. It can be checked that in any standard model
$\Mm(\F)$, if $\oo,\e,\xi$ are in the above described configuration,
then the ``standard" Minkowski-distances between $\oo,\e$ and
$\oo,\xi$ are the same. Conversely, if these two distances agree
then there exists an event $\e'$ as in Figure~\ref{med-fig}. Let now
$\e_1,\e_2$ be any two time-like separated events, $\e_1\prec\e_2$.
Then the relativistic distance between $\e_1,\e_2$ is the same as
that between the ``parallel translations" $\oo,\e$ of these, where
the ``parallel translation" happens according to
Figure~\ref{med-fig} (where for $\e''$ it is important only that it
is connected to both $\oo$ and to $\e_1$ with a light-signal, e.g.,
it is not important that $\oo\prec\e''$). If $\e_1\succ\e_2$ then we
define $\mu_{\ea,\so}(\e_1,\e_2)=-\mu_{\ea, \so}(\e_2,\e_1)$.  We
note that, while this relativistic distance strongly depends on the
parameters $\ea,\oo$, the relativistic {\it equi\,}distance relation
we get from this does not depend on $\ea,\oo$ any more. So, let us
define relativistic equidistance, or 4-equidistance, as
\begin{description}
\item[]
$\Med(\e_1,\e_2,\e_3,\e_4)\deiff
\mu_{\ea,\so}(\e_1,\e_2)=\mu_{\ea,\so}(\e_3,\e_4)$,\quad for any
$\ea\in\Par$ and event $\oo$ on $\ea$'s world-line.
\end{description}
\medskip

Having defined the desired coordinate systems in $\SigTh$, we
conclude this section with some remarks on what this method can give
us, what it can be used for.

We asked earlier, in section~\ref{sec:sigth}, where the paradigmatic
effects --- a-synchronicity, time-dilation, length-contraction ---
of special relativity theory came into the picture in Signalling
theory. One answer is the following. We defined natural coordinate
systems to the particles. (These coordinate systems correspond to
the observers in $\SpecRel$, this correspondence will be made
explicit in section~\ref{sec:reduc}.) Now, the
coordinate-transformations between these are so that the three
paradigmatic effects of special relativity (mentioned in
section~\ref{sec:sr}) hold in a version where we can recalibrate the
units of measurement.

This section contains definitions only, definitions (with some
parameters) in the language of $\SigTh$ that in the standard models
define coordinate systems for the particles/experimenters. We can
get an axiom system characterizing the standard models (thus doing
the job of $\SigTh$) via using these definitions. Namely, we can
state as axioms that the coordinate systems defined for the
experimenters have all the good properties we want (e.g., the
beginning and end-points of light-signals are exactly those of the
ordered segments of slope 1). This alternative axiom system would be
more complicated and less natural than $\SigTh$ of \cite{Ax},
however, it would be the result of a clear-cut method that can be
used in many other situations, as indicated below.

We can use the method of this section for exploring space-times
other than the special relativistic one, and for using signals of
various different nature, too. We mention some examples briefly, we
think that elaborating these examples would be worthwhile.

We can use the method of setting up a coordinate system as described
in this section, for example, for a particle moving
faster-than-light (FTL) in a special relativistic space-time. So,
let us take as standard models the standard models $\Mm(\F)$
modified so that the particles are the lines with slope more than 1
(and not the ones with slope less than 1). If we apply our method to
these modified models, then the FTL experimenter $\ea$ will find
that its space $\Space_{\ea}$ is a 3-dimensional Minkowski-space
$\MS\de\langle \F^3,\Bw,\Med\rangle$, and not a Euclidean space
$\langle \F^3,\Bw,\Ed\rangle$. He can reach by signals directly, and
check whether they are motionless w.r.t.\ him, only those
places/brothers that are time-like separated from him in terms of
$\MS$, but he can get indirect information about the rest of places
by communicating with these primarily reachable brothers. By working
through the details, we can get an axiom system
$\SigTh^{\mbox{\small ftl}}$ axiomatizing the signalling models of
FTL experimenters that would be quite analogous to Ax's $\SigTh$.
The main difference would be that the first group of axioms for
3-dimensional Euclidean space would be replaced by an analogous
axiom system for 3-dimensional Minkowski space. For this we can use
the one devised by Goldblatt in \cite[Appendix A]{Gold}. For a
slightly different approach for including FTL observers in this
setting see \cite{hoffman}.

However, communicating with directed signals (as in $\SigTh$)
between FTL experimenters is rather restricted if we want to take
the experiments to be executable (FTL experimenters can get
information this way only about the part of their space $\MS$ which
is in their ``past" in terms of $\MS$ as a Minkowski-space). We can
change the nature of signals to be undirected (but otherwise letting
their speed to be 1), imagining that if two events are connected
with a signal, then the information this signal carries appears at
both events ``at once". This is connected somehow to time-travel, a
subject strongly connected to FTL motion. The method given in the
present section is suitable for exploring space-time with undirected
signals, too.

The method given in this section can also be used for giving meaning
to two-dimensional time. Time being 2-dimensional could simply mean
that the events happening with the experimenters can be best
described by, say, the structure $\langle\F^2,\prec\rangle$. For
example, one could assume that our experimenter lives in a world
characterized by the 2+2-dimensional Minkowski-metric
$\sqrt{t_1^2+t_2^2-x^2-y^2}$ and then apply our method to see what
kind of coordinate system he would set up for himself, and in
general, what kind of responses he would get to his experiments.

Finally, we can imagine using signals of infinite velocity, this way
we can explore the Newtonian space-time characterized by absolute
time. Or, we can use bent signals of general relativity. For
example, we can explore the outer part of the Schwarzshild black
hole (the space-time outside the event horizon) with the same
method. We would take as experimenters a team of densely placed
suspended observers (spaceships in outer space using their drives to
maintain their desired positions), constantly checking positions by
communicating with photons (as light-signals), and using
freely-falling spaceships (or astronauts) as messengers.

\section{Defining new entities, interpretations}\label{sec:defth}
Our aim is to clarify the connections between $\SpecRel$ and
$\SigTh$. Not only the vocabularies of these two theories are
disjoint, even on the intuitive level they speak of different kinds
of things. We can see that somehow photons and observers of
$\SpecRel$ correspond to signals and particles of $\SigTh$, but what
correspond to quantities in $\SigTh$? Quantities of $\SpecRel$ do
not seem to enter the picture in $\SigTh$. Yet, in section
\ref{sec:algo} we defined something that intuitively could
correspond to quantities in $\SigTh$. In this section we recall some
tools from mathematical definability theory by which we can make
explicit the way quantities arise in $\SigTh$.

We briefly recall the tools that we will use in the next section for
making connections between theories for special relativity in a very
precise sense. We elaborated these tools in \cite{AMNDef,Mad02} for
the specific purpose of establishing a strong connection between two
versions of special relativity theory, the so-called
observer-independent geometrical and the reference-frame oriented
ones. We only recall the syntactic form to be used in specifying a
concrete interpretation together with some background intuition. We
elaborated a more extensive definability theory for this kind of
connecting theories that we do not recall here. We will say some
words about it at the end of this section. For simplicity, we will
treat function symbols as special relation symbols.

In ``traditional", one-sorted definability theory, an interpretation
of a theory $\Th'$ in language $\LL'$ into a theory $\Th$ in another
language $\LL$ is the following. For each $n$-place relation symbol
$R$ of $\LL'$ we assign a formula $\varphi_R$ of $\LL$ with at most
$n$ free variables. (We think of $\varphi_R$ as the ``definition of
$R$\," in $\LL$.) This then defines a natural translation function
$\tr:\LL'\longrightarrow\LL$ by replacing each atomic formula
$R(v_1,\dots,v_n)$ with $\varphi_R(v_1,\dots,v_n)$. This is an
interpretation of $\LL'$ into $\LL$. This interpretation is an
interpretation of $\Th'$ into $\Th$ iff $\Th$ proves the translated
theory $\Th'$, i.e.,
\begin{description}
\item[$\star$]
$\Th\models\tr(\psi)\mbox{\  whenever\  }\Th'\models\psi,\quad\mbox{
for all }\psi\in\LL'$.
\end{description}

\noindent On the semantic side, an interpretation of $\Th'$ into
$\Th$ ``constructs" a model of $\Th'$ inside each model of $\Th$.
Namely, it associates a model $\tr(\Mm)$ of $\LL'$ to each model
$\Mm$ of $\LL$ in such a way that the universe of $\tr(\Mm)$ is the
same as that of $\Mm$, and for each assignment $k$ of the variables
into this universe we have
\begin{description}
\item[$\star\star$]
$\tr(\Mm)\models\psi[k]\quad\mbox{ if and only if
}\quad\Mm\models\tr(\psi)[k]$,\quad for each formula $\psi$ in
$\LL'$.
\end{description}

In the new, ``non-traditional" or ``generalized" definability theory
we will use a notion of interpretation that does the same thing,
except that the universe of $\tr(\Mm)$ will not necessarily be a
subset of the universe of $\Mm$, therefore its definition and the
property analogous to ($\star\star$) above will be more involved. We
will define new entities as elements of new ``sorts". Using
many-sorted FOL is not an essential feature of this generalized
definability theory, just it is convenient in many cases, as it is
in our present task.

We illustrate the idea of defining new sorts with a simple example.
The language of affine planes in, e.g., \cite{Gold} is two-sorted,
we have two sorts $\Points, \Lines$ and we have a binary relation
between them, the relation $I$ of incidence (or membership) between
a point and a line. Another language in use for the same is
one-sorted, see, e.g., \cite{TGgeom}, we have one sort $\Points$ and
we have a three-place relation $\Col$ of ``collinearity" between
three points. Everyone can connect the two ways of thinking about
affine planes immediately: a line is the set of all points collinear
with given two distinct points. Thus a line $\ell$ is a subset of
the old universe, given two distinct points $p,q$ the line $\ell$
going through them is defined by
\begin{description}
\item[]
$\ell(p,q)\de\{ x : \Col(x,p,q)\}$.
\end{description}
But the new sort $\Lines$ stands for the set of all these subsets!
We can specify one line with the open formula $\Col(x,p,q)$ with one
free variable $x$, but how can we define the set of all lines? Well,
we will define the set of the parameters $p,q$ specifying the
individual lines: we identify the set of all lines with the set of
pairs of distinct points. Thus the formula defining the new sort
$\Lines$ will have two free variables $p,q$ and it will state $p\ne
q$. We are almost there, except that different pairs of distinct
points may specify the same line, and we have to take this into
account when talking about equality of lines, i.e., when
interpreting the equality symbol on the sort $\Lines$. We can do
this again with a formula using 4 free variables $p,q,p',q'$ stating
when the lines specified by $p,q$ and $p',q'$ coincide. In our case
this formula can be taken to be $\Col(p',p,q)\land\Col(q',p,q)$.

So far we have defined the universe of the new sort $\Lines$ and the
equality relation of this new sort by two formulas in the ``old"
language, i.e., in the language talking about $\Points$ and $\Col$.
Having defined a universe means that we have variables ranging over
this universe (and we can quantify over them). In other words, we
have to introduce variables $\Var(\Lines)$ of sort $\Lines$. Then,
in order to be able to use the definition of the new sort $\Lines$,
we need to connect $\Var(\Lines)$ to variables used in the
definition for $\Lines$, i.e., to $\Var(\Points)$. We can state this
connection by matching a variable $\ell$ of sort $\Lines$ with
variables denoting its ``defining parameters", e.g., we can state
that $\ell_p,\ell_q$ denote parameters that define $\ell$. After
this we can define the incidence relation, too: $I(x,\ell)\deiff
\Col(x,\ell_p,\ell_q)$, where $x$ is a variable of sort $\Points$
and $\ell$ is a variable of sort $\Lines$.

Summing up: defining the new sort $\Lines$ goes by defining the
variables $\Var(\Lines)$ of the new sort and matching them to the
variables of the old sort $\Var(\Points)$ occurring in the defining
formula of the sort $\Lines$, defining the equality on the sort
$\Lines$, and defining the non-logical symbol of incidence $I$ which
involves the sort $\Lines$. Thus we can interpret the 2-sorted
language of affine planes into the one-sorted one by the following
data:
\begin{description}
\item[]
$\var: \ell \mapsto \langle\ell_p,\ell_q\rangle$\quad for
$\ell\in\Var(\Lines)$,
\item[]
$\Lines(\ell)\deiff \ell_p\ne\ell_q$,
\item[]
$\ell=h\deiff \Col(\ell_p,h_p,h_q),\Col(\ell_q,h_p,h_q)$,
\item[]
$I(x,\ell)\deiff \Col(x,\ell_p,\ell_q)$.
\end{description}

\noindent The above data then define a translation function $\tr$
from the 2-sorted language of affine planes to their one-sorted
language as follows:

\begin{description}
\item[]
$\tr(\exists \ell\psi)\de \exists \ell_p,\ell_q\enskip\ell_p\ne\ell_q,\tr(\psi)$,\\
$\tr(\ell=h)\de \Col(\ell_p,h_p,h_q),\Col(\ell_q,h_p,h_q)$,\\
$\tr(I(x,\ell))\de\Col(x,\ell_p,\ell_q)$,\\
 the rest of the definition of $\tr$ is more or less straightforward.
\end{description}

\noindent The new feature in this translation function, over the
traditional one, is that we translate the quantifiers according to
the defining formula and variable-matching of the new sort and we
translate equality on the new sort, too. Throughout, we will use the
above variable matching $\var: \ell \mapsto
\langle\ell_p,\ell_q\rangle$ without recalling it.

This translation is not only recursive and structural, it is also
meaning preserving in the sense analogous to ($\star\star$). In more
detail: let $\Mm=\langle P,\Col\rangle$ be a model of the one-sorted
language. We will construct its ``translation", a model $\tr(\Mm)$
of the two-sorted language. Let
\begin{description}
\item[]
$U\de \{ \langle x,y\rangle\in P\times P : x\ne y \}$, and let
$E\subseteq U\times U$ be defined by
\item[]
$E\de\{ \langle u,v\rangle\in U\times U: \Col(u_1,v_1,v_2),
\Col(u_2,v_1,v_2)\}$.
\end{description}
Assume that $E$ is an equivalence relation on $U$, then define
\begin{description}
\item[]
$\tr(\Mm)\de\langle P,L,I\rangle$\quad where\quad\\
$L\de U/E$\quad and\\
$I\de \{ \langle x,u\rangle\in P\times L :  \Col(x,v_1,v_2) \mbox{
for some }v\in u/E\}$.
\end{description}
Let $\Var_P,\Var_L$ denote the sets of variables in the 2-sorted
language of the affine planes and let $\Var_P'\de
\Var_P\cup(\Var_L\times\{1\})\cup(\Var_L\times\{2\})$ be the
variables of the one-sorted language.  Now, let
$k:\Var_P\cup\Var_L\longrightarrow\tr(\Mm)$ be any evaluation of the
variables of the 2-sorted language, and let
$\tr(k):\Var_P'\longrightarrow\Mm$ be an  evaluation of the
variables of the one-sorted language such that $\tr(k)(x)=k(x)$ if
$x\in\Var_P$, and if $\ell\in\Var_L$ then $\langle
\tr(k)(\ell,1),\tr(k)(\ell,2)\rangle$ is an arbitrary element of
$k(\ell)$. Then the following is true for each formula $\psi$ of the
2-sorted language:

\begin{description}
\item[($\star\star'$)]
$\tr(\Mm)\models\psi[k]\quad\mbox{ if and only if }\quad
\Mm\models\tr(\psi)[\tr(k)]$.
\end{description}

\noindent The above ($\star\star'$) expresses that the translation
function preserves meaning when we talk about the 2-sorted model
constructed inside the one-sorted model.

Now, such a translation $\tr$ is an interpretation from $\Th'$ into
$\Th$ iff, just as before,
\begin{description}
\item[($\star'$)]
$\Th\models\tr(\psi)\mbox{\  whenever\  }\Th'\models\psi,\quad\mbox{
for all }\psi\in\LL'$.
\end{description}

Definitional equivalence of theories $\Th',\Th$ in different
languages $\LL',\LL$ is a strong connection between them, much
stronger than mutual interpretability requiring that the two
interpretations be inverses of each other, up to isomorphism. (Cf.\
\cite[Ex.4.3.46, p.266]{Mad02}.)

Two theories $\Th'$ and $\Th$ are said to be {\it definitionally
equivalent} if they have a common definitional extension. Here, two
theories are said to be the same if they prove the same formulas.
But what is a definitional extension? In the one-sorted case,
definitional extension of $\Th$ is $\Th\cup\Delta$ where $\Delta$ is
a union of definitions of the form $\Delta(R)\de
\{R(v_1,\dots,v_n)\leftrightarrow\f_R(v_1,\dots,v_n)\}$ with $\f_R$
as above ($\star$) (see, e.g., \cite[pp.60-61]{Hodges}). For telling
what definitional extension is in the many-sorted case, we return to
our previous example of defining the sort $\Lines$. Let us write
$\delta(p,q)$ and $\e(p,q,p',q')$ for $p\ne q$ and
$\Col(p',p,q),\Col(q',p,q)$ respectively, for the formulas defining
the ``domain" and the ``equality" on the new sort $\Lines$. The
explicit definition of the sort $\Lines$ will also involve a new
relation $\pi$ fixing the connection of the new sort to the old
ones. Now, $\Delta(\Lines,\pi)$ is defined to be the set of the
following sentences

\begin{description}
\item[]
$\exists p,q\,(\pi(p,q,\ell),\pi(p,q,\ell'))\leftrightarrow
\ell=\ell'$,
\item[]
$\exists\ell\,(\pi(p,q,\ell),\pi(p',q',\ell))\leftrightarrow
\e(p,q,p',q')$,
\item[]
$\exists\ell\,(\pi(p,q,\ell))\leftrightarrow\delta(p,q)$.
\end{description}

\noindent We note that the intuitive meaning of $\pi(p,q,\ell)$ is
that ``$p,q$ are distinct points lying on $\ell$\,", or, ``$p,q$
code, or represent, line $\ell$\,". So far it was the variable
matching that played this role and, intuitively, $\pi(p',q',\ell)$
is an explicit way of saying $\e(p',q',\ell_p,\ell_q)$.

After having defined the new sort $\Lines$, the definition
$\Delta(I)$ of the incidence relation is the same as in the
one-sorted case:

\begin{description}
\item[]
$I(p,\ell)\leftrightarrow \exists
p',q'\,(\pi(p',q',\ell),\Col(p,p',q'))$.
\end{description}

\noindent Now,  $\Th\cup\Delta(\Lines,\pi)\cup\Delta(I)$ is a
definitional extension of $\Th$, where $\Th$ is the ``one-sorted"
theory of affine planes. A {\it definitional extension} of any
theory $\Th$ is  $\Th\cup\Delta$ where $\Delta$ is a union of
definitions of the above form. Instead of describing the above in
more detail, we refer to \cite{AMNDef}, \cite[sec.4.3]{Mad02},
\cite[sec.6.3]{pezs} where many examples can also be found.

The notion of definitional equivalence is important for our
purposes, and we believe that it is an important one in
understanding how we form our concepts. We try to illustrate this
with an example. We will see that the theory $\EFd$ of Euclidean
fields and the theory $\SigTh$ of special relativity are mutually
interpretable into each other. However, they are not definitionally
equivalent.%
\footnote{Similar observations apply to a slight variant
$\SpecRel_0+\Compl$ of $\SpecRel$ in place of $\SigTh$ (cf.\
Thm.\ref{mainthm} in section \ref{sec:defeq}). This can be extended
to the Newtonian theory in \cite[sec.4.1, p.423]{pezs}.} Namely,
$\SigTh$ and $\EFd$ cannot have a common definitional extension
because of the following two reasons. (i) $\SigTh$ has to be an
``information-losing" reduct of any definitional extension of
$\EFd$, and (ii) any theory is an ``information-preserving" reduct
of any of its own definitional extensions. We note that (ii) holds
because the very idea of ``definitional extension" is an extension
based on ``information" contained in the unextended theory; thus by
forgetting this extra structure we lose nothing, we can recover it
from the unextended theory. We explain (i): In a definitional
extension of $\EFd$ of which $\SigTh$ is a reduct, we will define
the new sort $\Par$ of experimenters together with a projection
function $\pi_P$ which ties the behavior of $\Par$ to $\EFd$. Such a
projection function will single out the experimenter whose
world-line is the time-axis, in other words, we can single out
``the" motionless experimenter.%
\footnote{The easiest way of making this precise is that there are
fields with no automorphisms at all, e.g., the field of real
numbers, and this means that the structure
$\langle\Par,\pi_P,\EFd\rangle$ will have no automorphism, either.}
Absolute motion! However, the {\it essence} of relativity theory is
that motion is relative. This is formalized in the so called Special
Principle of Relativity, which states that all the experimenters are
equivalent, we cannot tell which one is motionless and which one
moves. Indeed, any experimenter can be taken to any other
experimenter by an automorphism, in any model of $\SigTh$. Thus,
when making the reduct of a definitional extension of $\EFd$ in
order to obtain $\SigTh$, we have to forget $\pi_P$, otherwise we do
not get the right concept of experimenter. This is
``information-loss" since we cannot recover $\pi_P$ from $\SigTh$.
This shows that ``forgetting" is an important part in forming the
concept of experimenter in this case. ``Less is more" in this case.
Definitional equivalence keeps track of these kinds of
``forgetting", while mutual interpretability may not do this.
\bigskip

We conclude this section with a few words about interpretations. We
already wrote about the philosophical importance of interpretations
between theories in the introduction. Here we write about more
technical aspects. An interpretation $\tr$ from theory $\Th'$ to
$\Th$ is a connection between them, and this connection imports some
properties of one theory to the other. For example, if $\Th$ is
consistent, then $\Th'$ is also consistent. If $\tr$ is faithful and
$\Th'$ is undecidable, then $\Th$ is also undecidable, and if $\Th'$
and $\Th$ are mutually interpretable in each other, then an axiom
system for $\Th$ can be imported to $\Th'$ via any two mutual
interpretations. Definitional equivalence induces a strong duality
between $\Th'$ and $\Th$. For these kinds of application of
interpretations see, e.g., \cite{HFried, Tarski53, Pam04, Mad02}.
The present paper intends to show the usefulness of interpretations
in physical theories, e.g., defining operational semantics for a
physical theory. We note that definability theory is quite
extensively used in geometry, see, e.g., \cite[App.B]{Gold},
\cite{Pam04, Pam05, Pam07, TGgeom}.

Versions of the general interpretability we use in this paper
appeared in various different forms as early as in 1969, see
\cite{Pre69, Szcze77, MakRey77, BentP84, Mak85, Mak93, Hodges}.
Almost all of these works use a syntactic device similar to ours,
let's call it explicit definitions, but they all elaborate on
different semantical aspects of this general definability. For
example, \cite{BentP84} characterizes when a functor of a given form
is the semantical part of an interpretation. \cite{MakRey77, Mak85}
recast model theory in a categorical form, where both the
syntactical and semantical parts of an interpretation are functors
between pretoposes, and it is proved that both functors are
equivalences when one is. This theorem is called a conceptual
completeness theorem. For the model theoretical forms, meaning and
impacts of this completeness theorem we refer to \cite{harnik}. (We
refer specifically to \cite[section 6, item (3)]{harnik} for
connections with the notion of general interpretability.) In
\cite{AMNDef,Mad02}, it is shown that our form of explicit
definitions outlined in this section is not ad hoc in the sense that
any sensible definition can be brought to this form. Namely, a
notion of implicit definability suggests itself as a necessary
condition for these new entities to be called ``defined", see, e.g.,
\cite[sec.4.3.1]{Mad02} and \cite{Hodges}. An analogue of the Beth
definability theorem (\cite[4.3.48]{Mad02}) states that if a sort of
new elements is implicitly definable, then it is explicitly
definable, too. We note that the powerset of the universe of an
infinite model is not implicitly definable in the sense of
\cite[sec.4.3.1]{Mad02}, while, say, the set of two-element subsets
of it is implicitly (and thus also explicitly) definable.

We hope that the content of this section is enough to give us a
guiding intuition for what comes in the rest of this paper.

\section{Reducing SpecRel to Signalling theory: an interpretation}
\label{sec:reduc} In this section we define in detail an
interpretation of $\SpecRel_0$ in $\SigTh$. We have to define (over
$\SigTh$) the new sorts $\Q$ and $\B$, and the new operations and
relations $+,\star,\Obs,\Ph,\W$ that involve these new sorts.

We begin with defining the new sort $\Q$. In section~\ref{sec:algo}
we already defined a field $\F(\ea,\oo,\uu)$, that will provide the
definition to our new sort $\Q$ and to $+,\star$. However, that
definition had three parameters $\ea,\oo,\uu$ (the particle who was
setting up his coordinate system, the ``beginning of the era", and
duration of one year). Up to isomorphism, we get the same field no
matter how we choose these 3 parameters, but their universes
strongly depend on the parameters (namely, the universe of
$\F(\ea,\oo,\uu)$ is the set of events on $\ea$'s world-line). Which
one should we take as the set of elements of sort $\Q$? The answer
is: take neither one, take all of them! Intuitively, this means that
we take the disjoint union of all the fields belonging to the
different parameters, and then we define an equivalence relation on
this set that relates the isomorphic images of the same element. For
this, in our explicit definition of the quantity sort we need a
uniform formula that defines the isomorphisms between the fields
$\F(\ea,\oo,\uu)$. One such formula is given in \cite[p.305]{Mad02}.
Here we give a simpler formula defining the isomorphisms between the
various incarnations of our field. We can give this simpler formula
because relativistic equidistance is available for us, while
\cite{Mad02} used only the betweenness relation.

\begin{figure}[h!]
\begin{center}
\psfrag{a}[b][b]{$\ea$} \psfrag{a'}[b][b]{$\ea'$}
\psfrag{o}[t][l]{$\oo$} \psfrag{u}[b][r]{$\uu$}
\psfrag{x}[b][r]{$\x$} \psfrag{o'}[l][l]{$\oo'$}
\psfrag{u'}[l][l]{$\uu'$} \psfrag{x'}[l][l]{$\x'= \x''\slash \uu''$}
\psfrag{x''}[l][l]{$\x''=\mu_{\ea',\so'}(\oo,\x)$}
\psfrag{u''}[l][l]{$\uu''=\mu_{\ea',\so'}(\oo,\uu)$}
\includegraphics[scale=.3]{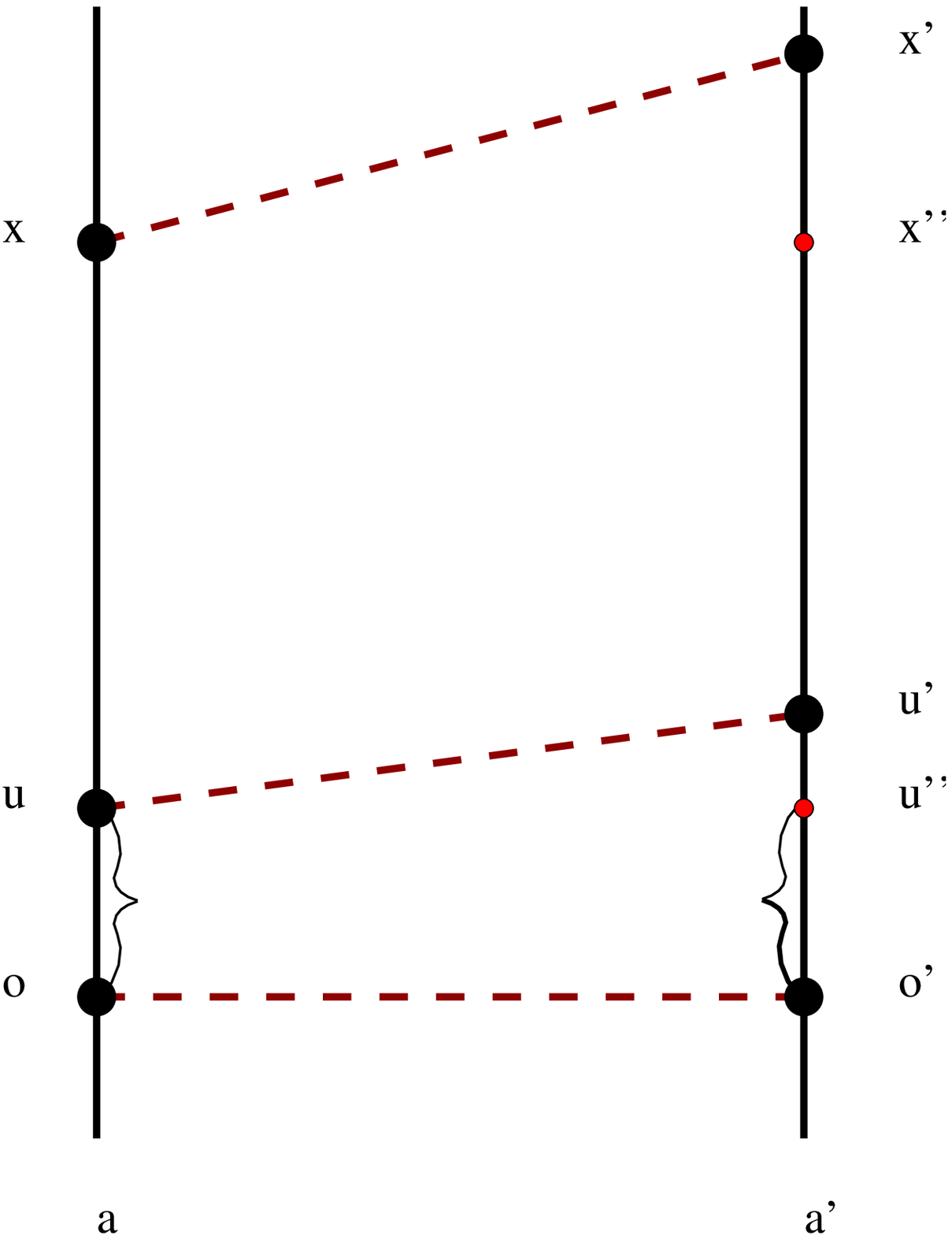}
\end{center}
\caption{\label{iso-fig} The isomorphism
$\f_{\is}(\ea,\oo,\uu,\ea',\oo',\uu')$ between $\F(\ea,\oo,\uu)$ and
$\F(\ea',\oo',\uu')$.}
\end{figure}

We are going to define the isomorphisms sought for between the
fields $\F(\ea,\oo,\uu)$. See Figure~\ref{iso-fig}. Let
$\ea,\oo,\uu,\ea',\oo',\uu'$ be suitable parameters for defining the
fields (as in section~\ref{sec:algo}). The isomorphism between them
will take $\oo$ to $\oo'$, $\uu$ to $\uu'$ and it will take an
arbitrary $\x$ on the world-line of $\ea$ to $\x'\de \x''\slash
\uu''$ where $\x'',\uu''$ are events on $\ea'$'s world-line such
that $\Med(\x,\oo,\x'',\oo')$ and $\Med(\uu,\oo,\uu'',\oo')$,
further $\slash$ denotes the division operation of the field
belonging to $\ea',\oo',\uu'$. Let
$\f_{\is}(\x,\x',\ea,\oo,\uu,\ea',\oo',\uu')$ denote the formula
expressing the above. We denote the isomorphism as
$\f_{\is}(\ea,\oo,\uu,\ea',\oo',\uu')$, and we denote the unique
$\x$ with the property $\f_{\is}(\x,\x',\ea,\oo,\uu,\ea',\oo',\uu')$
as $\x=\f_{\is}(\x',\ea,\oo,\uu,\ea',\oo',\uu')$.

Let $\Fp(\ea,\oo,\uu)$ express that $\ea,\oo,\uu$ are appropriate
parameters for a field $\F(\ea,\oo,\uu)$, let $U$ be the disjoint
union of the universes of all the fields $\F(\ea,\oo,\uu)$, and let
$E$ denote the binary relation relating isomorphic elements, i.e.,

\begin{description}
\item[]
$\Fp(\ea,\oo,\uu)\deiff \Ev(\oo), \Ev(\uu), \oo\ne\uu, \oo\prec\uu,
\ea\T\oo, \ea\T\uu$,
\item[]
$U\de\{ \langle\x,\ea,\oo,\uu\rangle : \Fp(\ea,\oo,\uu),
\x\in\F(\ea,\oo,\uu)\}$,
\item[]
$E\de\{\langle(\x,\ea,\oo,\uu),(\x',\ea',\oo',\uu')\rangle :
\f_{\is}(\x,\x',\ea,\oo,\uu,\ea',\oo',\uu')\}$.
\end{description}

It can be shown that $E$ is an equivalence relation on $U$, in each
standard model of $\SigTh$. Our quantity sort will be $U/E$.

Recall that we are in the process of defining $\SpecRel_0$ over
$\SigTh$.

We are ready to define the quantity sort $\Q$ explicitly, by using
the tools we introduced in the previous section. If $q$ is a
variable of the (new) sort $\Q$, then
$q_{\x},q_{\ea},q_{\oo},q_{\uu}$ denote the corresponding variables
of the (old) sorts $\Sig$ and $\Par$. We can think of this variable
matching as $q$ denotes an equivalence block of $E$ (i.e., an
element of $U/E$), and $\langle
q_{\x},q_{\ea},q_{\oo},q_{\uu}\rangle$ denotes an arbitrary
(unknown) element in the equivalence block $q$. Intuitively, $q$
denotes an ``abstract" quantity, and
$\f_{\is}(q_{\x},q_{\ea},q_{\oo},q_{\uu},\ea,\oo,\uu)$ is the
corresponding ``concrete" quantity in the field $\F(\ea,\oo,\uu)$.
Let us denote this last thing as
\begin{description}
\item[]
$\rep(q,\ea,\oo,\uu)\de\f_{\is}(q_{\x},\ea,\oo,\uu,q_{\ea},q_{\oo},q_{\uu})$.
\end{description}
This situation is somewhat analogous to the concept of a manifold in
general relativity theory, the elements of the manifold are the
``observer-independent" entities, and the charts/observers associate
concrete values to these. Below comes the definition of the sort
$\Q$:

\begin{description}
\item[]
$\var: q\mapsto\langle q_{\x},q_{\ea},q_{\oo},q_{\uu}\rangle$\quad
for $q\in\Var_{\Q}$.
\item[]
$\Q(q)\deiff q_{\ea}\T q_{\x}, \Ev(q_{\x}),
\Fp(q_{\ea},q_{\oo},q_{\uu})$,
\item[]
$q=q'\deiff
\f_{\is}(q_{\x},q_{\x}',q_{\ea},q_{\oo},q_{\uu},q_{\ea}',q_{\oo}',q_{\uu}')$.
\end{description}

\noindent Note that this definition of the sort $\Q$ is analogous to
the one given for the new sort $\Lines$ in the example of affine
planes in the previous section.

We get the definitions for $+,\star$ from writing up the definitions
given in section~\ref{sec:algo}, as follows.  Recall the formula
$+(\st,\st_1,\st_2,\ea,\oo)$ from section~\ref{sec:algo}.

Now, here is the definition of addition of sort $\Q$:

\begin{description}
\item[]
$+(q,q_1,q_2)\de \\
+(q_{\x},\f_{\is}(q_{1\x},q_{\ea},q_{\oo},q_{\uu},q_{1\ea},q_{1\oo},q_{1\uu}),
\f_{\is}(q_{2\x},q_{\ea},q_{\oo},q_{\uu},q_{2\ea},q_{2\oo},q_{2\uu}),q_{\ea},q_{\oo})$.
\end{description}

\noindent The formula defining multiplication of sort $\Q$ is
obtained analogously.

The rest of this section (interpreting $\SpecRel_0$ in $\SigTh$)
will be relatively straightforward.

We turn to defining the sort $\B$. We will define the sort $\B$ of
bodies as the union of observers and photons. So first we define the
entities that we will call photons. A photon will be defined just as
a signal $\s$ that is not an event. The world-line of this photon
will be defined as the set of all events that lie on the
4-dimensional line defined by the beginning and end points of $\s$.
This way, many photons will share the same world-line, just as in
the case of affine planes many pairs of distinct points define the
same line, and we will define two photons to be equal if they share
the same world-line. An observer will be defined to be a coordinate
system. We recall from section~\ref{sec:algo} that six parameters
are required for defining a coordinate system, namely the
experimenter $\ea$, a ``zero" $\oo$ and a time-unit $\uu$, and three
locations $a_x,a_y,a_z$ specifying the space coordinate axes. These
parameters have to satisfy the conditions below, which we will
denote by $\Op$ ($\Op$ refers to ``observer parameters"):

\begin{description}
\item[]
$\Op(\ea,\oo,\uu,a_x,a_y,a_z)\deiff \Fp(\ea,\oo,\uu),  \ea\| a_x,
\ea\| a_y, \ea\| a_z, \Ort(\ea,a_x,\ea,a_y), \Ort(\ea,a_x,\ea,a_z),
\Ort(\ea,a_y,\ea,a_z)$.
\end{description}

\noindent Two observers will be defined equal if they assign the
same coordinates to all events.

We are ready to formalize these definitions by using the tools we
introduced in section~\ref{sec:defth}. Let $\Var_{\B}$ denote the
set of variables of sort $\B$. If $b$ is a variable of sort $\B$,
then $b_{\s},b_{\ea},b_{\oo},b_{\uu},b_{x},b_{y},b_{z}$ will denote
the corresponding variables of ``old" sorts. Intuitively, this body
will be $b_{\s}$ if this is a ``real", non-degenerate signal (i.e.,
if $b_{\s}$ is not an event), and if $b_{\s}$ is ``degenerate"
(i.e., if it is an event), then the body $b$ will be the observer
$\langle b_{\ea},b_{\oo},b_{\uu},b_{x},b_{y},b_{z}\rangle$. We are
ready to define the new sort $\B$ together with the unary formulas
$\Ph(b)$ and $\Obs(b)$:

\begin{description}
\item[]
$\var: b\mapsto\langle
b_{\s},b_{\ea},b_{\oo},b_{\uu},b_{x},b_{y},b_{z}\rangle$\quad for
$b\in\Var_{\B}$.
\item[]
$\Ph(b)\deiff \lnot\Ev(b_{\s})$,
\item[]
$\Obs(b)\deiff
\Ev(b_{\s}),\Op(b_{\ea},b_{\oo},b_{\uu},b_{x},b_{y},b_{z})$,
\item[]
$\B(b)\deiff\Ph(b)\lor\Obs(b)$,
\end{description}

\noindent We are going now to define the equality relation on this
new sort $\B$. For stating equality of photons, first we express
that three events are on one light-like line
($\lambda(\e_1,\e_2,\e_3)$), then we express that an event is on the
world-line of a signal ($\wl(\e,\s)$).

\begin{description}
\item[]
$\lambda(\e_1,\e_2,\e_3)\deiff
\bigwedge\{\exists\s[(\e_i,\s,\e_j)\lor(\e_j,\s,\e_i)] :
i,j\in\{1,2,3\}\}$,
\item[]
$\wl(\e,\s)\deiff\exists \e_1,\e_2\enskip
\lambda(\e,\e_1,\e_2),\Beg(\s,\e_1),\End(\s,\e_2)$.
\end{description}

\noindent Recall from section~\ref{sec:algo} that the formula
$\cord(\e,\ea,\oo,\uu,a_x,a_y,a_z)=(\st,\g_x, \g_y, \g_z)$ expresses
that the coordinates of the event $\e$ are $\st, \g_x, \g_y, \g_z$,
in the coordinate system specified by $\ea,\oo,\uu,a_x,a_y,a_z$ .

\begin{description}
\item[]
$b=b'\deiff \\
(\lnot\Ev(b_{\s}),\lnot\Ev(b'_{\s}),\forall \e\,
\wl(\e,b_{\s})\leftrightarrow\wl(\e,b'_{\s}))\lor \\
(\Ev(b_{\s}),\Ev(b'_{\s}),\forall \e\,
\cord(\e,b_{\ea},b_{\oo},b_{\uu},b_x,b_y,b_z)=
\cord(\e,b'_{\ea},b'_{\oo},b'_{\uu},b'_x,b'_y,b'_z))$.
\end{description}

It remains to define the world-view relation $\W$. The intuitive
meaning of the formula $\W(m,b,t,x,y,z)$ will be that $m$ is an
observer, and the event at place $t,x,y,z$ in $m$'s coordinate
system is on the world-line of $b$. Let $m,b$ be variables of sort
$\B$ and let $t,x,y,z$ be variables of sort $\Q$. Assume that $m$ is
an observer, i.e., $\Ev(m_{\s})$. Let us denote the concrete value
of an abstract quantity $q$ in $m$'s coordinate system by
\begin{description}
\item[]
$m(q)\de\rep(q,m_{\ea},m_{\oo},m_{\uu})$.
\end{description}

\noindent We can now define $\W$ as follows:

\begin{description}
\item[]
$\W(m,b,t,x,y,z)\deiff  \exists\e\,
\cord(\e,m(t),m(x),m(y),m(z),m_{\ea},m_{\oo},m_{\uu},m_x,m_y,m_z),\\
(\lnot\Ev(b_{\s})\to\wl(\e,b_{\s})), (\Ev(b_{\s})\to b_{\ea}\T\e),
\Ev(m_{\s})$.
\end{description}

By the above, we gave definitions for all the sort and relation
symbols of the language of $\SpecRel$ in the language of $\SigTh$.
This defines a translation function $\tr$ between the two languages.
Let $=_{\Q}$ and $=_{\B} $ stand for the equality relations between
terms of sort $\Q$ and $\B$, respectively. In the next theorem we
state, without proof, that we indeed obtained an interpretation.

\begin{theorem}\label{interpret}
$\tr$ as given in this section is an interpretation of $\SpecRel_0$
into $\SigTh$, that is, the following are true:
\begin{description}
\item[]
$\SigTh\models``=_{\Q}\mbox{ and }=_{\B}\mbox{ are equivalence
relations}"$,
\item[]
$\SigTh\models``\mbox{the formulas defining
}+,\star,\Ph,\Obs,\W\mbox{ are invariant under }=_{\Q},=_{\B}"$,
\item[]
$\SigTh\models\tr(\psi)$\quad for all $\psi\in\SpecRel_0$.\quad
\end{description}
\end{theorem}

Having defined the desired interpretation of $\SpecRel_0$ into
$\SigTh$, in the next section we extend this interpretation to a
definitional equivalence between a slightly stronger version of
$\SpecRel_0$ and $\SigTh$.

\section{Definitional equivalence between SpecRel and Signalling theory} \label{sec:defeq}

In this section we investigate interpretability and definitional
equivalence between some of the FOL theories formalizing special
relativity. We show that a slightly reinforced version of
$\SpecRel_0$ is definitionally equivalent to $\SigTh$. We mean
interpretability and definitional equivalence in the sense of the
generalized definability theory of \cite{AMNDef,pezs,Mad02} outlined
in section~\ref{sec:defth}.

The interpreted theory $\tr(\SpecRel_0)$ is stronger than the
original one in the sense that there are sentences $\psi$ in the
language of $\SpecRel_0$ such that $\SigTh\models\tr(\psi)$ while
$\SpecRel_0\not\models\psi$. Such a sentence is, e.g., ``all lines
of slope less than 1 are world-lines of observers".
 We can express exactly how much more is true in the translated models
 by amending $\SpecRel_0$ with some existence,
extensionality, and time-orientation axioms (see below) and showing
that the so obtained theory is definitionally equivalent with
$\SigTh$. This is what we are going to do now.

The formulas describing the ``difference" between $\SpecRel_0$ and
$\SigTh$ are as follows. Formulas expressing that we have all kinds
of possible observers (from each point, in each direction, for each
velocity less than the speed of light there is an observer moving in
that direction with that speed, each observer can re-coordinatize
its coordinate-system with any space-isometry, each observer can set
the unit of its clock arbitrarily), and otherwise we are as economic
as possible (at most one photon through any two distinct events,
only one observer with the same coordinate-system, only photons and
observers as bodies, only one time-orientation for each observer).

These additional axioms, except the one about setting the clocks,
are denoted as {\sf AxThEx}, {\sf AxCoord}, {\sf AxExtOb}, {\sf
AxExtPh}, {\sf AxNobody}, \AxUp\ in \cite[sec.\ 2.5]{AMNHbSL}. Let
{\sf AxClock} formulate that each observer can set the unit of its
clock arbitrarily (in the spirit of the above axioms). Let $\Compl$
denote the set of these axioms and let $\SpecRel_0^+$ denote the
theory $\SpecRel_0$ amended with these formulas:

\begin{description}
\item[]
$\Compl \de \{{\sf AxThEx}, {\sf AxCoord}, {\sf AxClock}, {\sf
AxExtOb}, {\sf AxExtPh}, {\sf AxNobody}, \AxUp\}$,
\item[]
$\SpecRel_0^+\de \SpecRel_0 + \Compl$.
\end{description}

To state definitional equivalence between $\SpecRel_0^+$ and
$\SigTh$, we now define an interpretation $\Tr$ of $\SigTh$ into
$\SpecRel_0^+$. We have to define the universes $\Par,\Sig$ of
particles and signals and the relations $\T,\R$ of transmitting and
receiving, inside $\SpecRel_0$. Intuitively, particles are defined
to be observers, with two particles being equal if their world-lines
coincide:

\begin{description}
\item[]
$\var: a\mapsto a_b$\quad for $a\in\Var_{\Par}$, where
$a_b\in\Var_{\B}$.
\item[]
$\Par(a)\deiff \Obs(a_b)$,
\item[]
$a=a'\deiff \forall t,x,y,z\enskip\W(a,a',t,x,y,z)\leftrightarrow
x=y=z=0$.
\end{description}

Signals are defined to be photons with two events on their
world-lines representing the beginning and end-points of the signal.
We represent the two events with observers meeting the photon. The
following formulae express in $\SpecRel_0$ that ``in $b$'s
world-view, $p$ meets $a$ at time $t$", and ``$a,p,e$ meet in one
event", respectively:

\begin{description}
\item[]
$\Meet(b,p,a,t)\deiff \exists
x,y,z\enskip\W(b,p,t,x,y,z),\W(b,a,t,x,y,z)$,
\item[]
$\meet(a,p,e)\deiff \exists
b,t\enskip\Meet(b,a,p,t),\Meet(b,a,e,t)$.
\end{description}

\noindent  Now we are ready to interpret signals in $\SpecRel_0$:

\begin{description}
\item[]
$\var: \s\mapsto \langle\s_b,\s_p,\s_e\rangle$\quad for
$\s\in\Var_{\Sig}$, where $\s_b,\s_p,\s_e\in\Var_{\B}$.
\item[]
$\Sig(\s)\deiff \Ph(\s_p),\Obs(\s_b),\Obs(\s_e),\exists t\le
t'\,\Meet(\s_b,\s_p,\s_b,t),\Meet(\s_b,\s_p,\s_e,t')$.
\item[]
$\s=\s'\deiff \meet(\s_b,\s'_b,\s_p), \meet(\s_e,\s'_e,\s_p),
\lnot\meet(\s_b,\s_p,\s_e)\to \s_p=\s'_p$.
\end{description}

\noindent Finally,

\begin{description}
\item[]
$a\T\s\deiff \meet(a_b,\s_b,\s_p)$,
\item[]
$a\R\s\deiff \meet(a_b,\s_e,\s_p)$.
\end{description}

The above define a translation function $\Tr$ as indicated in
section~\ref{sec:defth}. We state in the next theorem, without
proof, that this $\Tr$ interprets $\SigTh$ in $\SpecRel_0^+$, and
moreover, together with the interpretation $\tr$ defined in the
previous section it forms a definitional equivalence between
$\SpecRel_0^+$ and $\SigTh$. This is the main theorem of this paper.

\begin{theorem}\label{mainthm}
$\SigTh$ is definitionally equivalent to $\SpecRel_0 + \Compl$, the
pair $\tr,\Tr$ of interpretations forms a definitional equivalence
between them.
\end{theorem}

We can read the above theorem as saying that what the theory
$\SigTh$ tells about special relativity is exactly what the theory
$\SpecRel_0 + \Compl$ says. Since no axiom in $\Compl$ follows from
$\SpecRel_0$, we can conclude that $\SigTh$ tells more than
$\SpecRel_0$, the amount of ``more" is exactly $\Compl$. However, we
did not include the axioms of $\Compl$ into $\SpecRel$, because we
do not need them in proving the main predictions of relativity
theory; we feel that they do not belong to the core of the physical
theory. Moreover, of the axioms of $\Compl$, we consider only
$\AxUp$ as having a physically (or even philosophically) relevant
content, namely it says that ``time is oriented".

On the other hand, we will see that $\SpecRel$ has a content that
$\SigTh$ does not say about special relativity theory. This is the
axiom $\AxSym$ of $\SpecRel$. So, what is the connection between
$\AxSym$ and $\SigTh$? Below we answer this question.

The interpretation $\tr$ we defined in the previous section does not
interpret $\SpecRel$ in $\SigTh$, because $\tr(\AxSym)$ does not
follow from $\SigTh$ (i.e., it is not true in the standard models
$\Mm(\F)$ of $\SigTh$).%
\footnote{We note that $\SpecRel$ can be interpreted in $\SigTh$ in
the way that we interpret $\SpecRel$ in the field $\Q,+,\star$.}
The reason for this is the following. $\AxSym$ states that any two
observers use the same units of measurement. We can express in the
language of $\SigTh$ that ``two observers use the same units of
measurement", and this defines an equivalence relation on the set of
all observers. For $\AxSym$ to be true, we should select any one of
the blocks of this equivalence relation (since $\AxSym$ states that
any two observers use the same units of measurements). But which one
should we select? The question might sound familiar. In the previous
analogous case (that concerned the various incarnations
$\F(\ea,\oo,\uu)$ of the field $\F$) we took all the classes ``up to
isomorphism". However, in the present case there are no definable
bijections between the blocks of this equivalence relation.

We can get around this problem by adding to the models $\Mm(\F)$ of
$\SigTh$ a ``unit of measurement". We can do this, e.g., the
following way. We add a new basic two-place relation symbol $\Tu$
(short for ``Time unit") of sort $\Sig$ to the language of
Signalling theory. In each standard model $\Mm(\F)$ we interpret
$\Tu$ as the set of pairs of events with Minkowski-distance 1. (We
note that these relations are not definable in $\Mm(\F)$ in the
language of Signalling theory.) Let us denote the so expanded
standard models by $\Mm(\F)^+$, and let $\SigTh^+$ denote an axiom
system for their theory (in the extended language). Now, the
interpretation we gave in this section can be extended to interpret
$\SpecRel$ in $\SigTh^+$. Moreover, it also can be made into a
definitional equivalence between a stronger version of $\SpecRel$
and $\SigTh^+$, that we obtain from $\SpecRel_0^+$ by exchanging
{\sf AxClock} with $\AxSym$:

\begin{description}
\item[]
$\Compl^- \de \{{\sf AxThEx}, {\sf AxCoord}, {\sf AxExtOb}, {\sf
AxExtPh}, {\sf AxNobody}, \AxUp\}$,
\end{description}

To our minds, the following theorem clarifies the connection between
$\AxSym$ and $\SigTh$. It says that the content of $\AxSym$ is to
set the time-unit: the difference between $\SigTh^+$ and $\SigTh$ is
that in $\SigTh^+$ we can express Minkowski-distance, while in
$\SigTh$ we have only Minkowski equidistance.
\goodbreak

\begin{theorem}\label{symthm} The following (i),(ii) hold:
\begin{description}
\item[{\rm (i)}] $\SigTh^+$ is definitionally equivalent to
$\SpecRel + \Compl^-$.
\item[{\rm (ii)}]
$\SigTh^+$ is not definitionally equivalent to $\SigTh$.
\end{description}
\end{theorem}

The proof of part (i) of the above theorem goes by extending the
interpretation $\Tr$ to $\SigTh^+$, this amounts to defining the new
relation $\Tu$ in $\SpecRel$; and also making some (minor) changes
in the definition of $\tr$. The proof of part (ii) of the above
theorem goes by showing that the automorphism groups of members of
$\SigTh$ and $\SigTh^+$ differ from each other, this technique is
elaborated in, e.g., \cite{Mad02,Hodges}.

We included $\AxSym$ into $\SpecRel$ as a tool for convenience, it
seems to carry no philosophical or physical importance. $\AxSym$ is
only a simplifying assumption.

Concerning some of the other theories for special relativity, we
mention that $\SigTh$ is definitionally equivalent to Goldblatt's
theory for special relativity in \cite[Appendix A]{Gold} amended
with time-orientation. I.e., the two theories are almost the same,
the only difference is that $\SigTh$ assumes time-orientation while
Goldblatt's theory does not. The proof of this last statement can be
put together from the definitions and ideas in
sections~\ref{sec:algo}, \ref{sec:reduc}. Also, (a slight variant
of) our $\SpecRel$ is definitionally equivalent to (a slight variant
of) Suppes's axiomatization of special relativity in \cite{Sup59,
Sup}.

\section{Conclusion} \label{sec:concl}
We intended to show in this paper some results the methods of
mathematical logic can provide for other branches of science, in
particular, for physics and the methodology of science. Using the
tools of definability theory of first-order logic, we compared in
detail two rather different axiom systems for special relativity
theory.  One of these, $\SpecRel$ of \cite{AMNSz}, is
coordinate-system-, or reference frame-oriented, while the other,
$\SigTh$ of \cite{Ax}, uses meager resources and talks about
particles emitting and absorbing signals. The two theories use
disjoint languages and talk about different kinds of entities. Yet,
a precise comparison was made possible by using mathematical logic,
and we obtained the following: $\SigTh$ can express and states
everything that $\SpecRel$ does, except for the relativistic
(Minkowski) distance between events (implied by $\AxSym$ in
$\SpecRel$), while in addition it states time-orientation for
space-time together with some auxiliary simplifying axioms
($\Compl$). Informally,

\begin{description}
\item[]
$\SigTh$ = $\SpecRel$ - relativistic distance + time-orientation +
auxiliaries,
\end{description}

\noindent and a little more formally

\begin{description}
\item[]
$\SigTh$ + $\AxSym$ = $\SpecRel$ + $\Compl^-$.
\end{description}

A byproduct of these investigations is a concrete operational
semantics for special relativity theory. We believe that
interpreting one theory in another is a flexible methodology for
connecting physical theories with each other as well as with the
``physical reality".


\begin{thebibliography}{99.}%


\bibitem{AMNDef} Andr\'eka, H., Madar\'asz, J. X., N\'emeti, I.:
Defining new universes in many-sorted logic. Preprint, Mathematical
Institute of the Hungarian Academy of Sciences, Budapest, 2001.
93pp.

\bibitem{pezs} Andr\'eka, H., Madar\'asz, J. X., N\'emeti, I.: On the logical
structure of relativity theories. Alfr\'ed R\'enyi Institute of
Mathematics, Hungar. Acad. Sci., Budapest, Research Report, July 5,
2002, with contributions from A. Andai, G. S\'agi, I. Sain and Cs.
T\H oke. \url{http://www.math-inst.hu/pub/algebraiclogic/
Contents.html}. 1312 pp.

\bibitem{AMNHbSL} Andr\'eka, H., Madar\'asz, J. X., N\'emeti, I.:
Logic of space-time and relativity theory. In: Handbook of Spatial
Logics. Eds: Aiello, M. Pratt-Hartmann, I., and van Benthem, J.
Springer Verlag, 2007. pp.607-711.

\bibitem{Vienna} Andr\'eka, H., Madar\'asz, J. X., N\'emeti, I.,
N\'emeti, P., Sz\'ekely, G.: Vienna Circle and Logical Analysis of
Relativity Theory. In: The Vienna Circle in Hungary (Der Wiener
Kreis in Ungarn). Eds: M{\'a}te, A., R{\'e}dei, M., Stadler, F.
Veroffentlichungen des Instituts Wiener Kreis, Band 16, Springer
Verlag, 2011, pp.247-268.

\bibitem{AMNSz} Andr\'eka, H., Madar\'asz, J. X., N\'emeti, I., Sz\'ekely,
G.: A logic road from special relativity to general relativity.
Synthese 186,3(2012)633-649.

\bibitem{Ax} Ax, J.: The elementary foundations of spacetime.
Foundations of Physics, 8,7/8 (1978), 507-546.

\bibitem{Sneed} Balzer, W., Moulines, U., Sneed, J. D.: An
architectonic for science. The structuralist program. D. Reidel
Publishing Company, Dordrecht, 1987.

\bibitem{vB82} van Benthem, J.: The logical study of science.
Synthese 51 (1982), 431-472.

\bibitem{vB12} van Benthem, J.: The logic of empirical theories
revisited. Synthese 186,2 (2012), 775-792.

\bibitem{BentP84} van Benthem, J., Pearce, D.: A mathematical
characterization of interpretation between theories. Studia Logica
43,3 (1984), 295-303.

\bibitem{BG} Burstall, R., Goguen, J.: Putting theories together
to make specifications. In: Proc. IJCAI'77 (Proceedings of the 5th
International Joint Conference on Artificial Intelligence, Vol 2,
pp.1045-1058. Morgan Kaufmann Puboishers Inc. San Francisco, CA,
USA, 1977.

\bibitem{Carnap}  Carnap, R.: Die Logische Aufbau der Welt. Felix Meiner Verlag,
Leipzig, 1928.

\bibitem{FOM} Friedman, H.: On foundational thinking 1. FOM (Foundations of Mathematics)
Posting, Archives www.cs.nyu.edu, January 20, 2004.

\bibitem{HFried} Friedman, H.: Interpretations of Set Theory in
Discrete Mathematics and Informal Thinking. Lectures 1-3. Nineteenth
Annual Tarski Lectures, Berkeley, 2007.
\url{http://www.math.osu.edu/~friedman.8/}

\bibitem{Frie} Friedman, M.: Foundations of space-time theories.
Relativistic physics and philosophy of science. Princeton University
Press, 1983.

\bibitem{Gard} G{\"a}rdenfors, P., Zenker, F.: Theory change as
dimensional change: conceptual spaces applied to the dynamics of
empirical theories. Synthese 190 (2013), 1039-1058.

\bibitem{Gold} Goldblatt, R.: Orthogonality and spacetime geometry.
Springer-Verlag, 1987.

\bibitem{harnik} Harnik, V.: Model theory vs. categorical logic: two
approaches to pretopos completion (a.k.a. $T^{eq}$). Centre de
Recherches Math\'ematiques CRM Proceedings and Lecture Notes Vol 53.
American Mathematical Society 2011. pp.79-106.

\bibitem{Hodges} Hodges, W.: Model Theory. Cambridge University
Press, 1993.

\bibitem{hoffman} Hoffman, B.: A logical treatment of special
relativity, with and without faster-than-light observers. BA Thesis,
Lewis and Clark College, Oregon, USA. 63pp. 2013.
arXiv:1306.6004[math.LO].

\bibitem{Volt} Konev, B., Lutz, C., Ponomaryov, D., Wolter, F.:
Decomposing description logic ontologies. In: Proceedings of 12th
Conf.\ on the Principles of Knowledge Representation and Reasoning,
Association for the Advancement of Artificial Intelligence, 2010.
pp.236-246.

\bibitem{Lutz} Lutz, C., Wolter, F.: Mathematical logic for life
science ontologies. In: Proceedings of WOLLIC-2009, Ono, H.,
Kanazawa, M., de Queiroz, R. (eds.), LNAI 5514, Springer, 2009.
pp.37-47.

\bibitem{Mad02} Madar\'asz, J. X.: Logic and relativity (in the
light of definability theory). PhD Dissertation, E\"otv\"os Lor\'and
University, 2002.
\url{http://www.math-inst.hu/pub/algebraic-logic/diszi.pdf}

\bibitem{MadSzek13} Madar\'asz, J. X., Sz\'ekely, G.: Special
relativity over the field of rational numbers. International Journal
of Theoretical Physics 52,5 (2013), 1706-1718.

\bibitem{Mak85} Makkai, M.: Ultraproducts and categorical logic.
In: Methods in Mathematical Logic, Springer LNM 1130 (1985),
222-309.

\bibitem{Mak93} Makkai, M.: Duality and definability in first order
logic. Memoirs of the American Mathematical Society No 503, 1993.

\bibitem{MakRey77} Makkai, M., Reyes, G.: First order categorical
logic. Lecture Notes in Mathematics 611. 1977.

\bibitem{Pam04} Pambuccian, V.: Elementary axiomatizations of projective
space and of its associated Grassman space. Note de Matematica 24,1
(2004/05), 129-141.

\bibitem{Pam05} Pambuccian, V.: Groups and Plane Geometry. Studia Logica
81 (2005), 387-398.

\bibitem{Pam07} Pambuccian, V.: Alexandrov-Zeeman type theorems
expressed in terms of definability. Aequationes Math. 74 (2007),
249-261.

\bibitem{Pre69} Previale, F.: Rappresentabilit$\grave{a}$ ed
equipollenza di teorie assiomatiche I. Ann. Scuola Norm. Sup. Pisa
23,3 (1969), 635-655.

\bibitem{Schelb} Schelb, U.: Characterizability of free motion in
special relativity. Foundations of Physics 30,6 (2000), 867-892.

\bibitem{Schutz} Schutz, J. W.: Independent axioms for Minkowski
space-time. Longman, 1997.

\bibitem{Szabo} Szab{\'o}, L. E.: Empirical foundation of space and
time. In: EPSA07: Launch of the European Philosophy of Science
Association, Su{\'a}rez, M., Dorato, MM., R{\'e}dei, M. (Eds.),
Springer, 2009.
\url{http://phil.elte.hu/leszabo/Preprints/LESzabo-madrid2007-preprint.pdf}

\bibitem{Szcze77} Szczerba, L. W.: Interpretability of elementary theories.
In: Butts, R. E., Hintikka, J. (eds) Logic, foundations of
mathematics and computability theory (Proc.\ Fifth Internat.\
Congr.\ of Logic, Methodology and Philos.\ of Sci., Univ.\ Western
Ontario, London, Ont., 1975), Part I, 129--145, Reidel, Dordrecht
(1977)

\bibitem{Sze} Sz\'ekely, G.: First-order logic investigation of relativity
theory with an emphasis on accelerated observers. PhD Dissertation,
E\"otv\"os Lorand University, Faculty of Sciences, Institute of
Mathematics, Budapest 2009. 150pp. ArXiv:1005.0973[gr-qc].

\bibitem{Sze10} Sz\'ekely, G.: A geometrical characterization of the
twin paradox and its variants. Studia Logica 95,1-2 (2010), 161-182.

\bibitem{Sup59} Suppes, P.: Axioms for relativistic kinematics with
or without parity. In L. Henkin, A. Taarski, and P. Suppes, editors,
Symposium on the Axiomatic Method with Special Reference to Physics,
North Holland, 1959. pp.297-307.

\bibitem{Sup} Suppes, P.: Some open problems in the philosophy of
space and time. Synthese 24 (1972), 298-316.

\bibitem{Tarski36} Tarski, A.: Der Wahrheitsbegriff in den formalisierten Sprachen.
Studia Philosophica 1 (1936).

\bibitem{Tarski} Tarski, A.: What is elementary geometry? In L.
Henkin, P. Suppes, and A. Tarski, editors, The axiomatic Method with
Special Reference to Geometry and Phsics, North-Holland, Amserdam
1959. pp.16-29.

\bibitem{TGgeom} Tarski, A., Givant, S. R.: Tarski's system of
geometry. The Bulletin of Symbolic Logic 5,2 (1999), 175-214.

\bibitem{Tarski53} Tarski, A., Mostowski, A., Robinson, R. M.: Undecidable theories.
North-Holland, Amsterdam, 1953.

\bibitem{Wheeler} Taylor, E. F., Wheeler, J. A.: Spacetime
Physics. Freeman and Co., San Francisco, 1963, 1966.

\end{thebibliography}
\end{document}